\documentclass[10pt,aps,longbibliography,onecolumn,superscriptaddress]{revtex4-2}

\usepackage[T1]{fontenc}
\usepackage[utf8]{inputenc}
\usepackage[english]{babel}

\usepackage{amsfonts,amsbsy,amssymb,amsmath,graphicx,float}

\graphicspath{{figures/}}

\usepackage[linktoc=all,colorlinks=true,citecolor=blue,linkcolor=blue]{hyperref}

\begin{document}
	
	\title{Lagrangian Descriptors \\ and the Action Integral of Classical Mechanics}
	
	\author{V\'ictor J. Garc\'ia-Garrido}
	\email{vjose.garcia@uah.es}
	\affiliation{Departamento de F\'isica y Matem\'aticas, Universidad de Alcal\'a, \\ Alcal\'a de Henares, 28871, Spain.}
	
	\author{Stephen Wiggins}
	\email{s.wiggins@bristol.ac.uk}
	\affiliation{School of Mathematics, University of Bristol, \\ Fry Building, Woodland Road, Bristol, BS8 1UG, United Kingdom.}
	
	\begin{abstract}
		In this paper we bring together the method of Lagrangian descriptors and the principle of least action, or more precisely, of stationary action, in both deterministic and stochastic settings. In particular, we show how the action can be used as a Lagrangian descriptor. This provides a direct connection between Lagrangian descriptors and Hamiltonian mechanics, and we illustrate this connection with benchmark examples.
	\end{abstract}
	
	\maketitle
	

	
	\section{Introduction}
	\label{sec:intro}
	
	In this paper we show how the method of Lagrangian descriptors (LDs) can be incorporated with the principle of least action (or more precisely, of stationary action) to provide a version of Lagrangian descriptors that is directly connected to Hamiltonian mechanics. Lagrangian descriptors are a trajectory diagnostic that facilitate the exploration of phase space structures, and global dynamics in general, in dynamical systems. This diagnostic tool was first introduced for the analysis of transport and mixing processes in the context of fluid dynamics \cite{madrid2009,mancho2013lagrangian,mendoza2010} in order to diagnose the flow structures governing transport processes of interest. Since that time, applications to a broader collection of applications have been carried out. For example, in the field of chemical reaction dynamics Lagrangian descriptors have been recognised as a valuable technique to compute the phase space structures that govern reactivity of the system, see e.g. \cite{craven2015lagrangian, craven2016deconstructing, craven2017lagrangian, revuelta2019unveiling, agaoglou2019,Main2017,Main2018, Main2019, Bartsch2016, GG2020b, naik2020, krajnak2019, gonzalez2020, GG2020a, katsanikas2020a}. 
	
	Lagrangian descriptors have a number of advantages over other trajectory diagnostics. They apply equally well to dynamical systems in closed form, to systems with arbitrary time dependence, as well as to dynamical systems described as data sets. In the domain of calculation for the Lagrangian descriptor, they can reveal all hyperbolic trajectories and their stable and unstable manifolds in a single calculation. This provides the input for the application of rigorous theorems such as the existence theorem for normally hyperbolic invariant manifolds (NHIMs) and their stable and unstable manifolds \cite{wiggins_normally_2014} and the Smale-Birkhoff homoclinic theorem for the existence of chaotic dynamics \cite{wiggins2003applied}. These features of Lagrangian descriptors are discussed in detail in \cite{ldbook2020}, and an open source software package to carry out the types of computations we describe here is provided in \cite{aguilar2021ldds}. However, to go into these issues in more detail is outside the scope of this paper, which is to demonstrate the connection between extremals of the action integral and Lagrangian descriptors, and to benchmark this connection on known closed-form examples where explicit calculations can be carried out.
	
	The principle of least (``extremal'') action is treated in most advanced books on classical mechanics; see, for example, \cite{arnol2013mathematical, goldstein2002classical, landau2013mechanics}. An intuitive and elementary discussion of the principle of least action is given by Richard Feynman in his famous lecture series, see \url{https://www.feynmanlectures.caltech.edu/II_19.html}. Extremals of the action integral are trajectories of the dynamical equations of motion. Extrema of the Lagrangian descriptor field typically correspond to stable and unstable manifolds of normally hyperbolic invariant manifolds, or invariant tori for long time averages of Lagrangian descriptors. Hence, we believe that a connection between these two trajectory diagnostics has the potential to yield greater insights into global dynamics and phase space structure. We note that use of the action as a Lagrangian descriptor was described in \cite{ldbook2020}, and has also been used in \cite{meseguer19, gonzalez2020, montoya2020phase}. See also \cite{montoya2021revealing}.
	
	This paper is outlined as follows. In Section \ref{sec:LDaction} we describe how the method of Lagrangian descriptors can be defined in terms of the action integral of classical mechanics. Next, in Sec. \ref{sec:saddle} we show by analytical means that this technique detects the stable and unstable manifolds for the linear saddle system with one degree of freedom (DoF). Section \ref{sec:harmosc} discusses, making use of the harmonic oscillator, how this approach can be applied to obtain Kolmogorov-Arnold-Moser (KAM) tori in dynamical systems and their frequencies, which characterize regular quasiperiodic motion. Section \ref{sec:2dof} focuses on the application of LDs to a model Hamiltonian system with 2 DoF from chemical reaction dynamics. In Sec. \ref{sec:stochsys} we illustrate how this trajectory-based diagnostic tool can be successfully applied to explore the phase space of random dynamical systems, and we do so by analyzing the stochastic dynamics of the Duffing oscillator. Finally, in Sec. \ref{sec:conc} we present the conclusions of this work.

	
	\section{Lagrangian Descriptors Based on the Classical Action}
	\label{sec:LDaction}
	
	Our setting will be an $n$ degrees-of-freedom Hamiltonian system whose Hamiltonian function is given by:
	\begin{equation}
		\mathcal{H}(\mathbf{q},\mathbf{p}) = T(\mathbf{p}) + V(\mathbf{q}) \;,
		\label{class_ham}
	\end{equation}
	where $T$ denotes the kinetic energy, which is a quadratic function of the momenta $\mathbf{p} = (p_1,\ldots,p_n) \in \mathbb{R}^n$ in the form:
	\begin{equation}
		T(p_1,\ldots,p_n) = \sum_{i=1}^n \dfrac{p_i^2}{2 \, m_i} \;,
	\end{equation}
	and the potential energy function $V(\mathbf{q})$ depends only on the configuration coordinates $\mathbf{q} = (q_1,\ldots,q_n) \in \mathbb{R}^n$. The corresponding Hamiltonian vector field is given by:
	\begin{equation}
		\begin{cases}
			\dot{q}_i = \dfrac{\partial \mathcal{H}}{\partial p_i} = \dfrac{p_i}{m_i} \\[.4cm]
			\dot{p}_i = -\dfrac{\partial \mathcal{H}}{\partial q_i} = -\dfrac{\partial V}{\partial q_i}
		\end{cases}
		,\quad i = 1, \ldots, n \;.
		\label{eq:hamvf}
	\end{equation}
	In order to study the phase space structures that characterize transport in the $2n-1$ energy hypersurface where motion occurs due to energy conservation, we can construct a Lagrangian descriptor by means of the classical action (Maupertuis action) $S_0$ calculated along a trajectory $\mathbf{x}(t) =  \left(\mathbf{q}(t),\mathbf{p}(t)\right)$ of Hamilton’s equations of motion in the phase space of the system:
	\begin{equation}
		S_0\left[\mathbf{x}(t)\right] = \int^{\mathbf{q}_1}_{\mathbf{q}_0} \mathbf{p} \cdot d\mathbf{q}  \, ,
		\label{action_def}
	\end{equation}
	where $\mathbf{q}_0 = \mathbf{q}(t_0)$ and $\mathbf{q}_1 = \mathbf{q}(t_1)$ are the initial and final configuration space points in the trajectory. The action $S_0$ defines a natural metric in the phase space which we use to construct a Lagrangian descriptor. Notice that we can rewrite the action integral as follows:
	\begin{equation}
		S_0\left[\mathbf{x}(t)\right] = \int^{\mathbf{q}_1}_{\mathbf{q}_0} \mathbf{p} \cdot d\mathbf{q} = \int^{t_1}_{t_0} \mathbf{p} \cdot \dfrac{d\mathbf{q}}{d t} \, dt = \int^{t_1}_{t_0}  2 T \, dt \, .
		\label{action_ke}
	\end{equation}
	This yields the integral along the phase space trajectory $\mathbf{x}(t)$ of twice the kinetic energy, which is a non-negative scalar function. 
	
	Given any initial condition $\mathbf{x}_0 = \mathbf{x}(t_0)$, we construct its trajectory by evolving the phase space point both forward and backward in time for the time periods $[t_0,t_0+\tau_f]$ and $[t_0-\tau_b,t_0]$, respectively. The action-based Lagrangian descriptor evaluated along the trajectory starting at $\mathbf{x}_0$ is defined by the expression:
	\begin{equation}
		\mathcal{S}(\mathbf{x}_{0},t_{0},\tau_{f},\tau_{b}) = \int^{ \mathbf{q}_{f}}_{ \mathbf{q}_{b}} \mathbf{p} \cdot d\mathbf{q} = \int^{t_{0}+\tau_{f}}_{t_{0}-\tau_{b}} 2T \; dt \;,
		\label{eq:LD_S}
	\end{equation}
	where $\mathbf{q}_f = \mathbf{q}(t_0 + \tau_{f})$ and $\mathbf{q}_b = \mathbf{q}(t_0 - \tau_{b})$ represent the configuration space endpoints of the resulting trajectory. The Lagrangian descriptor scalar field in Eq. \eqref{eq:LD_S} can be naturally split into two terms:
	\begin{equation}
		\mathcal{S}(\mathbf{x}_{0},t_{0},\tau_{f},\tau_{b})  = \mathcal{S}^{(f)}(\mathbf{x}_{0},t_{0},\tau_{f}) + \mathcal{S}^{(b)}(\mathbf{x}_{0},t_{0},\tau_{b}) \;,
	\end{equation}
	where $\mathcal{S}^{(f)}$ and $\mathcal{S}^{(b)}$ correspond, respectively, to the forward and backward in time contributions of the action integral to the total Lagrangian descriptor $\mathcal{S}$. The forward and backward components are given by the formulas:
	\begin{equation}
		\mathcal{S}^{(f)}(\mathbf{x}_{0},t_{0},\tau_{f}) = \int^{ \mathbf{q}_{f}}_{\mathbf{q}_{0}} \mathbf{p} \cdot d\mathbf{q} = \int^{t_{0}+\tau_{f}}_{t_{0}} 2T \; dt \quad,\quad \mathcal{S}^{(b)}(\mathbf{x}_{0},t_{0},\tau_{b}) =  \int^{\mathbf{q}_{0}}_{\mathbf{q}_{b}} \mathbf{p} \cdot d\mathbf{q} =  \int^{t_{0}}_{t_{0}-\tau_{b}} 2T \; dt \;.
		\label{eq:LD_S_fwbw}
	\end{equation}
	
	When the action-based Lagrangian descriptor function in Eq. \eqref{eq:LD_S} is calculated on a given grid of initial conditions in the phase space, it produces a scalar field that has the capability of highlighting the location and geometry of the invariant stable and unstable manifolds associated to the normally hyperbolic invariant manifolds in the phase space of the system \cite{wig2016,naik2019a}. These structures can be easily identified with features of the Lagrangian descriptor output where the scalar field displays abrupt changes in its values. It is important to remark that, if we consider only the contribution that comes from the forward integration of trajectories, $\mathcal{S}^{(f)}$, this will provide information about the stable manifolds, while the backward term $\mathcal{S}^{(b)}$ is used to determine the unstable manifolds. If both terms of the Lagrangian descriptor are added to obtain $\mathcal{S}$, the output can be used to visually locate the presence of NHIMs in the system at the intersection of the stable and unstable manifolds. Recall that, in the case of a 2 DoF Hamiltonian system, NHIMs are unstable periodic orbits (UPOs) whose stable and unstable manifolds have the topology of $S^{1} \times \mathbb{R}$ and are known in the literature as spherical cylinders.
	

\section{The 1 DoF Hamiltonian Saddle}
\label{sec:saddle}

In this section we prove that the action-based LD recovers the stable and unstable manifolds for the equilibrium point at the origin in the 1 DoF linear saddle dynamical system defined by the Hamiltonian function:
\begin{equation}
\mathcal{H}(q,p) = \dfrac{\lambda}{2} \left(p^2 - q^2\right) \, ,
\end{equation}
where $\lambda > 0$. The evolution of this dynamical system is determined by Hamilton's equations of motion:
\begin{equation}
\begin{cases}
\dot{q} = \dfrac{\partial \mathcal{H}}{\partial p} = \lambda \, p \\[.4cm]
\dot{p} = -\dfrac{\partial \mathcal{H}}{\partial q} = \lambda \, q 
\end{cases} \,.
\label{ds_saddle}
\end{equation}
We know that, given the initial condition $\mathbf{x}(t_0) = (q_0,p_0)$, the analytic solutions are:
\begin{equation}
q(t) = q_0 \cosh(\lambda t) + p_0 \sinh(\lambda t) \quad,\quad p(t) = p_0 \cosh(\lambda t) + q_0 \sinh(\lambda t) \;.
\end{equation}
This dynamical system has a hyperbolic equilibrium point at the origin with stable ($\mathcal{W}^{s}$) and unstable ($\mathcal{W}^{u}$) manifolds defined by the following sets:
\begin{equation}
\mathcal{W}^{s}(0,0) = \left\lbrace (q,p) \in \mathbb{R}^2 \; \Big| \; q = -p \right\rbrace \quad,\quad \mathcal{W}^{u}(0,0) = \left\lbrace (q,p) \in \mathbb{R}^2 \; \Big| \; q = p \right\rbrace \;.
\end{equation} 

We start by calculating the forward component of the action-based Lagrangian descriptor. Since the system is autonomous, we can choose without loss of generality $t_0 = 0$. If we take into account the definition of $\dot{q}$ given by Hamilton's equations, we can write:
\begin{equation}
\begin{split}
\mathcal{S}^{(f)}(\mathbf{x}_0,\tau) & = \int_{t_0}^{t_0 + \tau} p \, \dot{q} \, dt = \lambda \int_{0}^{\tau} p^2 \, dt = \lambda \int_{0}^{\tau} \left(p_0 \cosh(\lambda t) + q_0 \sinh(\lambda t)\right)^2 \, dt = \\[.3cm]
& = \lambda \left[q_0^2 \int_{0}^{\tau} \sinh^2(\lambda t) \, dt + p_0^2 \int_{0}^{\tau} \cosh^2(\lambda t) \, dt + 2 q_0 \, p_0 \int_{0}^{\tau} \sinh(\lambda t) \cosh(\lambda t) \, dt \right] = \\[.3cm]
& = \lambda \left[q_0^2 \int_{0}^{\tau} \sinh^2(\lambda t) \, dt + p_0^2 \int_{0}^{\tau} \cosh^2(\lambda t) \, dt + q_0 \, p_0 \int_{0}^{\tau} \sinh(2\lambda t)  \, dt \right] \\[.3cm]
& = \dfrac{\lambda}{2} \left(p_0^2 - q_0^2\right) \tau + \dfrac{1}{4}\left(q_0^2 + p_0^2\right) \sinh(2\lambda\tau) + \dfrac{1}{2} q_0 \, p_0 \left(\cosh(2\lambda \tau) - 1\right) \;.
\end{split}
\end{equation}
We have obtained thus the expression:
\begin{equation}
\mathcal{S}^{(f)}(\mathbf{x}_0,\tau) = \dfrac{\lambda}{2} \left(p_0^2 - q_0^2\right) \tau + \dfrac{1}{4}\left(q_0^2 + p_0^2\right) \sinh(2\lambda\tau) + \dfrac{1}{2} q_0 \, p_0 \left(\cosh(2\lambda \tau) - 1\right) \,.
\label{sfw_exp}
\end{equation}
We will show next that $\mathcal{S}^{(f)}$ attains a minimum at the stable manifold of the system. Consider $p_0$ as fixed and differentiate with respect to $q_0$ to yield:
\begin{equation}
\dfrac{\partial \mathcal{S}^{(f)}}{\partial q_0} (\mathbf{x}_0,\tau) = - \lambda \tau q_0 + \dfrac{1}{2} q_0 \sinh(2\lambda\tau) + \dfrac{1}{2} p_0 \left(\cosh(2\lambda \tau) - 1\right) \,.
\end{equation}
Imposing that the derivative has to vanish we get:
\begin{equation}
q_0 = - \dfrac{\cosh(2\lambda\tau) - 1}{\sinh(2\lambda\tau) - 2\lambda \tau} \, p_0 = - G\left(\tau;\lambda\right)  p_0 \,.
\end{equation}
Notice now that if $\tau$ is sufficiently large, it is straightforward to show that:
\begin{equation}
\lim_{\tau \to \infty} G\left(\tau;\lambda\right) = 1 \,,
\end{equation}
and therefore we recover the stable manifold $q_0 = - p_0$. To check that this is in fact a minimum we differentiate again:
\begin{equation}
\dfrac{\partial^{\, 2} \mathcal{S}^{(f)}}{\partial^2 q_0} (\mathbf{x}_0,\tau) = - \lambda \tau  + \dfrac{1}{2}\sinh(2\lambda\tau) \geq 0 \; , \quad \forall \; \tau > 0 \;.
\end{equation}
To check that this is indeed the case, it suffices to verify that the function is increasing in $\tau$, or equivalently that its derivative with respect to $\tau$ is non-negative:
\begin{equation}
\dfrac{\partial^{\, 3} \mathcal{S}^{(f)}}{\partial \tau \partial^2 q_0} = \lambda \left[\cosh(2\lambda\tau) - 1\right] = 2 \lambda \sinh^2(\lambda\tau) > 0 \;.
\end{equation}

At this point we calculate the value attained by the forward component of the action on an initial conditions lying on the stable manifold, that is $p_0 = - q_0$. This gives:
\begin{equation}
\mathcal{S}^{(f)}(q_0,-q_0,\tau) = \dfrac{1}{2} q_0^2 \left[\sinh(2\lambda\tau) - \cosh(2\lambda\tau)\right] + \dfrac{1}{2} q_0^2 \;.
\end{equation}
Taking the limit of this expression as $\tau$ gets very large gives:
\begin{equation}
\lim_{\tau \to \infty} \mathcal{S}^{(f)}(q_0,-q_0,\tau) = \dfrac{1}{2} q_0^2 \;.
\end{equation}
Therefore, the limit value that the forward component of LDs attains on the stable manifold as time goes to infinity is bounded for points on the stable manifold. On the other hand, if the initial condition $(q_0,p_0)$ is selected at a point that is not on the stable manifold, but close to it, the value attained by the forward component of LDs, given by Eq. \eqref{sfw_exp}, blows up when $\tau$ gets very large. This demonstrates that for sufficiently large values of the integration time $\tau$, the scalar field provided by $\mathcal{S}^{(f)}$ will display an abrupt change, or 'singular feature', at the points on the stable manifold. At those points, the LD function will become singular (or non-differentiable), and this property is relevant for the visualization of the manifolds when the LD scalar field is displayed, and also for the efficient extraction of these geometrical objects by means of image-processing techniques in the post-processing of numerical simulations.

We can go further with our analysis to show that a relationship can be obtained between the Lyapunov exponent $\lambda$ and the integration time $\tau$ for which the action-based LDs converges to the stable manifold. This would occur when:
\begin{equation}
G(\tau;\lambda) \approx 1 \quad \Leftrightarrow \quad \cosh(2\lambda\tau) \approx \sinh(2\lambda\tau) \quad \Leftrightarrow \quad \cosh(2\lambda\tau) - \sinh(2\lambda\tau) \approx 0  \;.
\end{equation}
Therefore, we can write:
\begin{equation}
\cosh(2\lambda\tau) - \sinh(2\lambda\tau) = e^{-2\lambda \tau} \approx 0 \quad \Leftrightarrow \quad e^{-2\lambda \tau} \approx 10^{-N} \quad \Leftrightarrow \quad \tau = \dfrac{N \ln(10)}{2\lambda} \,.
\end{equation}
Consequently, if we choose for example $N = 8$, which represents simple precision in the computer, we arrive at:
\begin{equation}
\tau = \dfrac{4 \ln(10)}{\lambda} \approx \dfrac{9.21}{\lambda} \;.
\end{equation}
In what follows, we will carry out all the simulations of LDs for the linear saddle system in Eq. \eqref{ds_saddle} using $\lambda = 1$ and integration times of $\tau = 6$ and $\tau = 8$, which would yield approximations to the stable and unstable manifolds of the order of $10^{-6}$.

\begin{figure}[htbp]
	\begin{center}
		A)\includegraphics[scale=0.24]{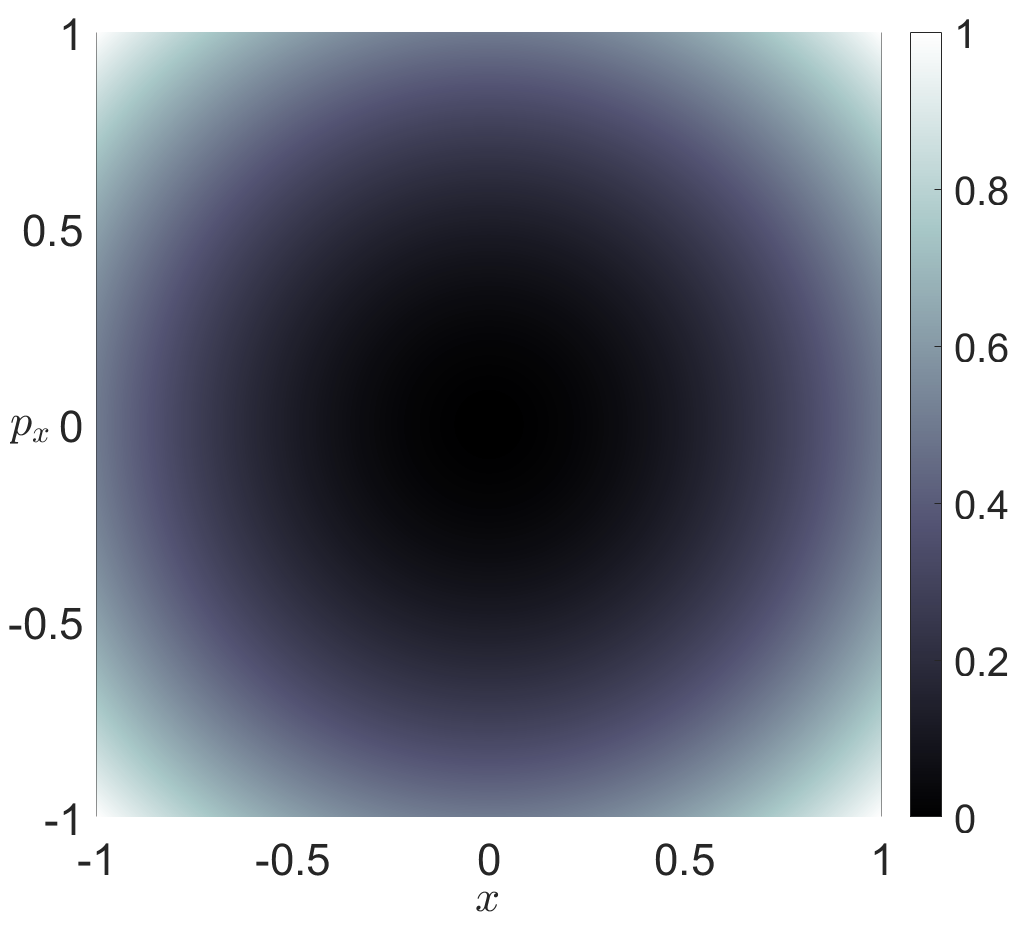}
		B)\includegraphics[scale=0.24]{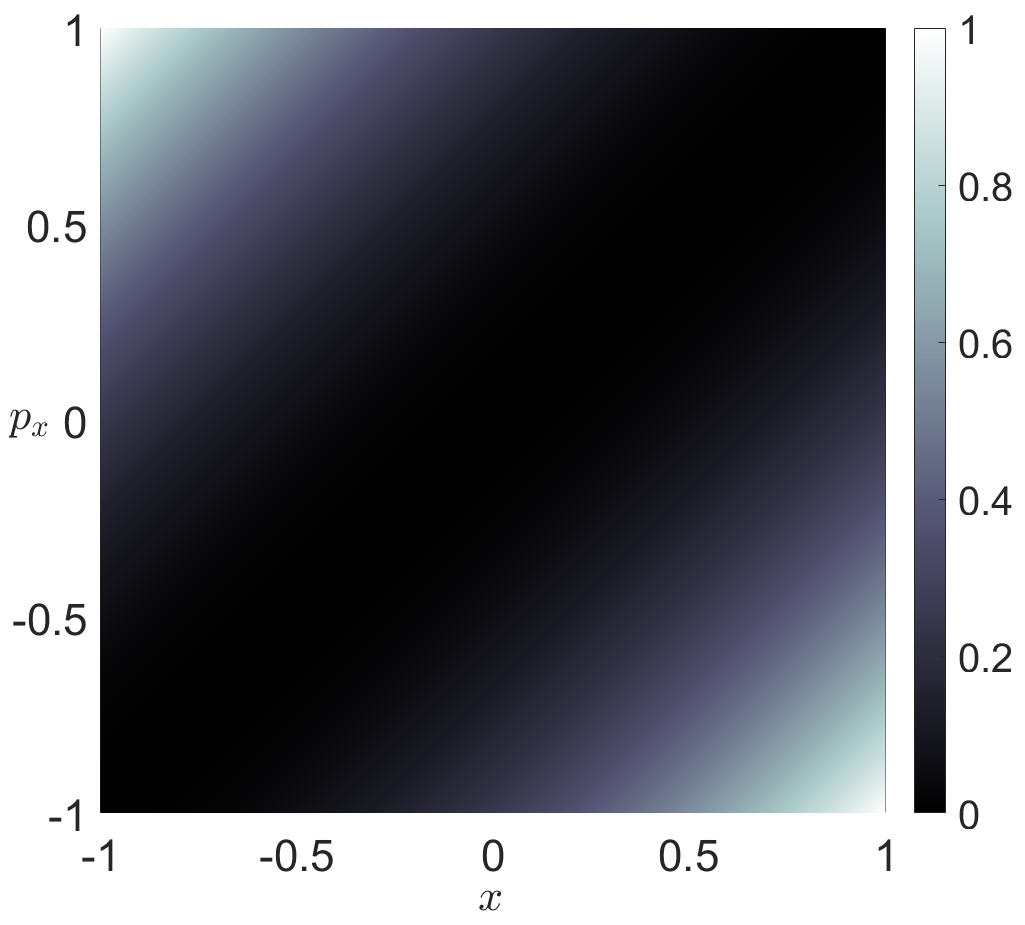}
		C)\includegraphics[scale=0.24]{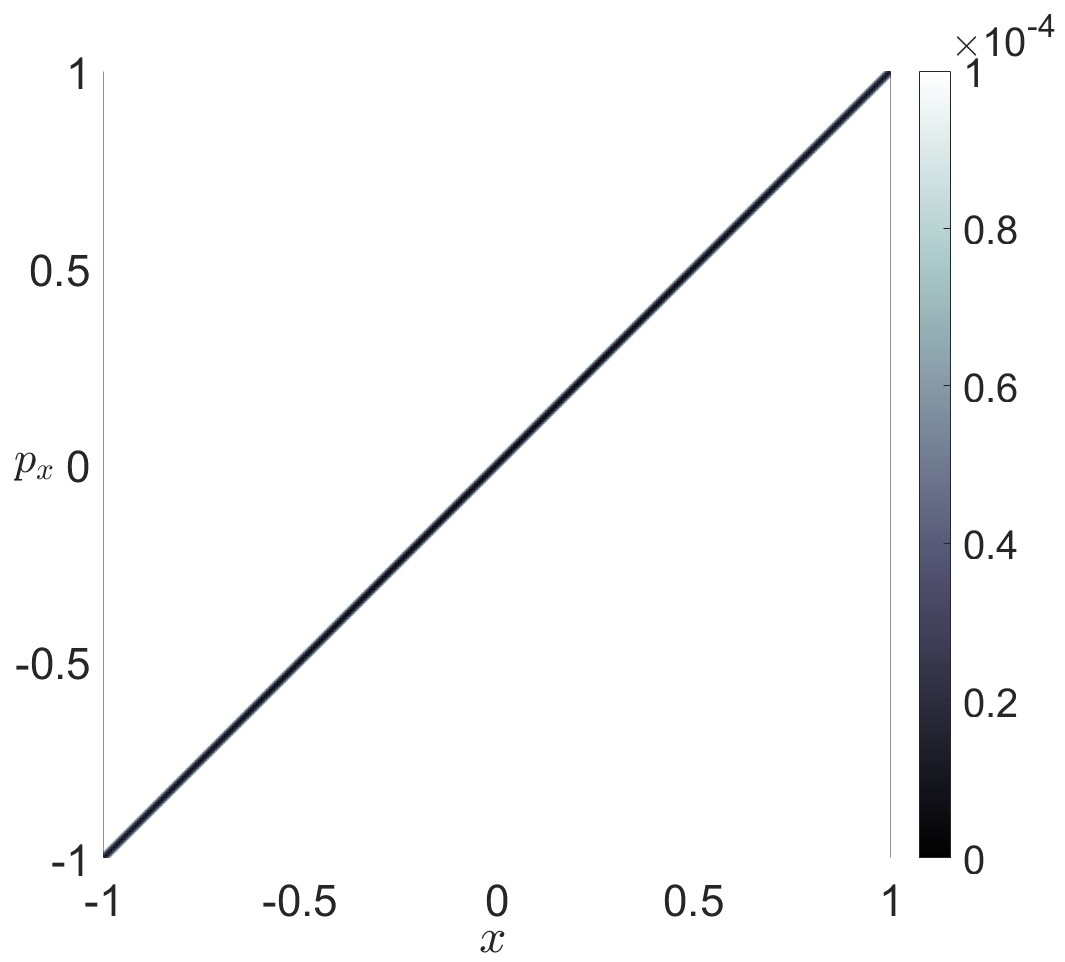}
		\end{center}
	\caption{A) Total LD with $\tau = 6$ for the linear saddle Hamiltonian system in Eq. \eqref{ds_saddle} using $\lambda = 1$. B) Backward LD for $\tau = 6$. For the computation of the LD scalar field, all the trajectories starting at the initial conditions on the grid have been integrated for the same time $\tau = 6$. C) The same as B) but the color scale has been adjusted to highlight the unstable manifold.}
	\label{ld_saddle_fixTime}
\end{figure}

It is important to remark here that all the arguments we have given above for $\mathcal{S}^{(f)}$ carry over naturally to the backward component of LDs, $\mathcal{S}^{(b)}$. It is a straightforward calculation to show that $\mathcal{S}^{(b)}$ is bounded at points on the unstable manifold of the system, and also that it attains a minimum at them. Moreover, the backward LD scalar field displays a sharp transition, i.e. a 'singular feature' at the unstable manifold. At this point, we would like to highlight that some difficulties arise when one tries to visualize both the stable and unstable manifolds simultaneously in the same plot for the 1 DoF linear saddle system. This becomes clear when the total LD scalar field is calculated, that is, the forward plus backward components are added up. If we do so, we get the following analytic expression:
\begin{equation}
    \mathcal{S}(\mathbf{x}_0,\tau) = \mathcal{S}^{(f)}(\mathbf{x}_0,\tau) + \mathcal{S}^{(b)}(\mathbf{x}_0,\tau) = \lambda \tau \left(p_0^2 - q_0^2\right) + \dfrac{\left(p_0^2 + q_0^2\right)}{4} \sinh(2\lambda \tau) = 2 \lambda \tau \mathcal{H}_0 + \dfrac{\left(p_0^2 + q_0^2\right)}{4} \sinh(2\lambda \tau)  \;,
\end{equation}
and the expression above evidences that for $\tau \gg 1$, the scalar field of the action-based LD behaves as:
\begin{equation}
    \mathcal{S}(\mathbf{x}_0,\tau) \approx \dfrac{\left(p_0^2 + q_0^2\right)}{4} \sinh(2\lambda \tau)  \;,
\end{equation}
which does not capture both manifolds simultaneously in the same plot. This is illustrated in Fig. \ref{ld_saddle_fixTime} A) for $\tau = 6$, where the total LD is displayed (the output has been normalized by subtracting the minimum value of the LD function, and then dividing by the difference between the maximum and minimum values of the LD). Notice however that if we plot for example the backward component of LDs, $\mathcal{S}^{(b)}$, see Fig. \ref{ld_saddle_fixTime} B) and C), the unstable manifold is visible. To enhance the detection of the unstable manifold we need to adjust the range of the color scale in the plot, since all initial conditions that are not on the unstable manifold yield very large values of $\mathcal{S}^{(b)}$, which makes it difficult to highlight in Fig. \ref{ld_saddle_fixTime} B) the location of the unstable manifold. This issue has been adequately fixed in Fig. \ref{ld_saddle_fixTime} C).

\begin{figure}[htbp]
	\begin{center}
		A)\includegraphics[scale=0.26]{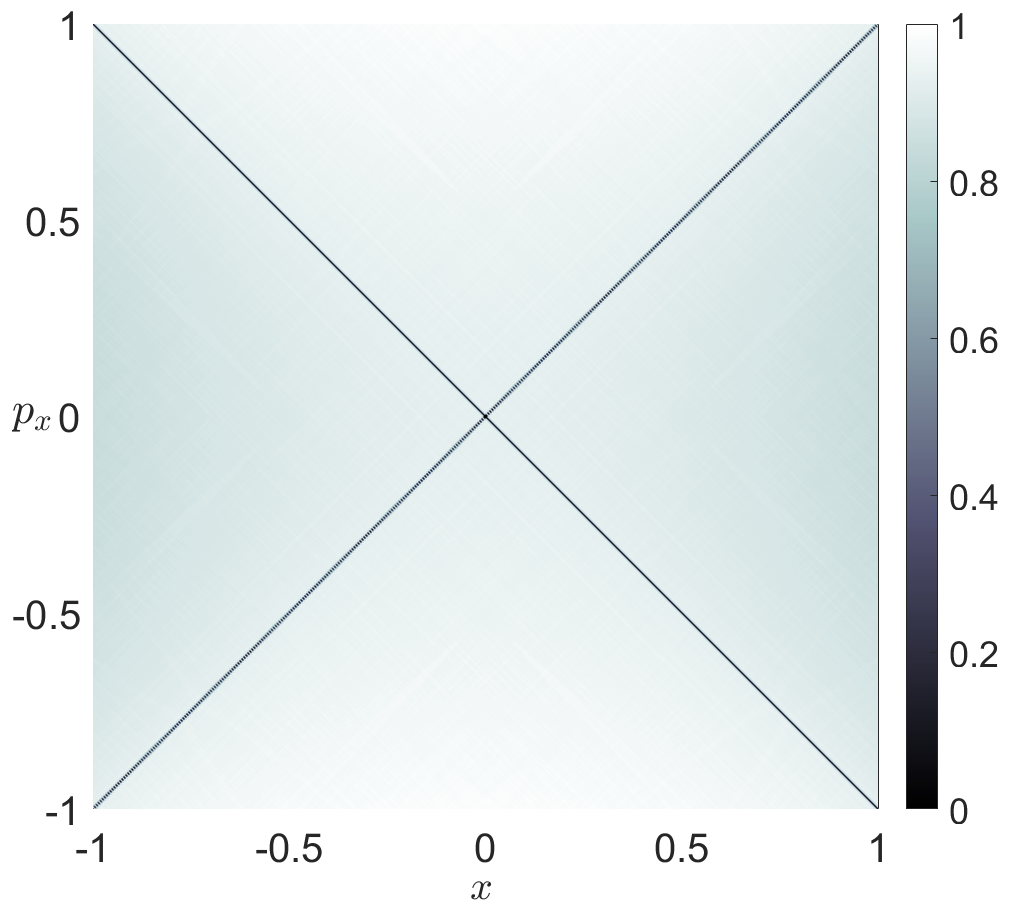}
		B)\includegraphics[scale=0.26]{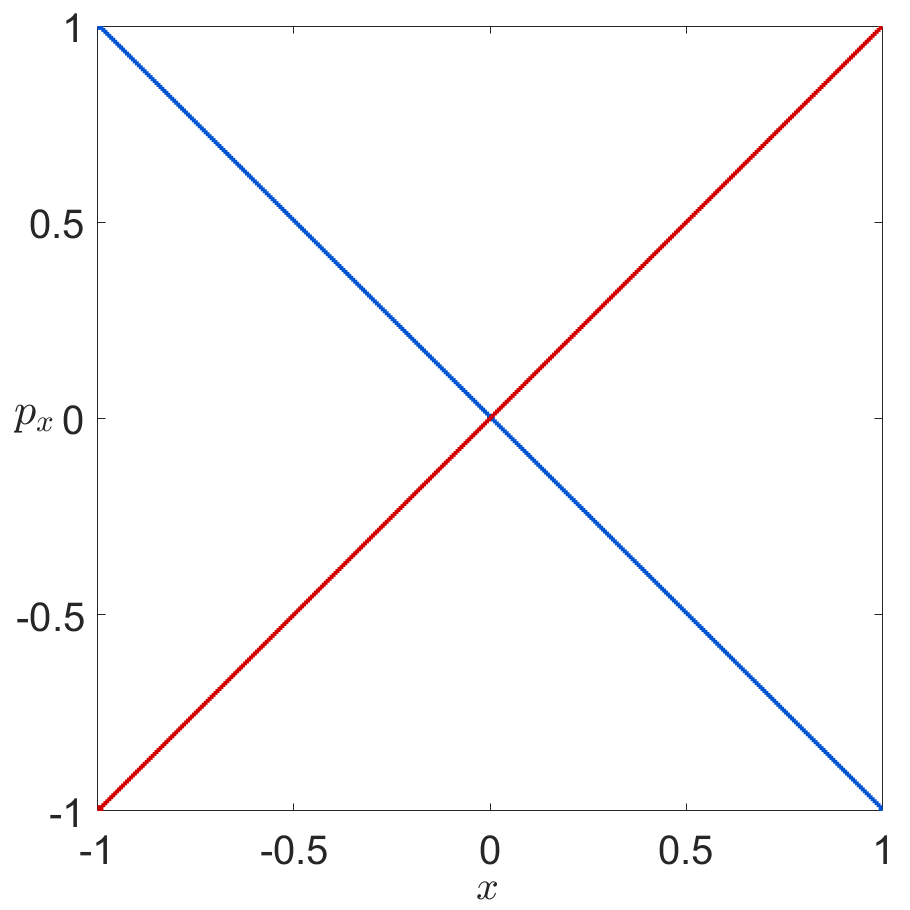}
		\end{center}
	\caption{A) Total LD with $\tau = 8$ for the linear saddle Hamiltonian system in Eq. \eqref{ds_saddle} using $\lambda = 1$. B) Stable (blue) and unstable (red) manifolds extracted from the LD scalar field displayed on the left panel. Trajectories have been integrated both forward and backward in time for $\tau = 8$ or until they leave a square $[-16,16]\times[-16,16]$ centered about the origin. This approach is known in the literature as variable time LD.}
	\label{ld_saddle_varTime}
\end{figure}

The reason why this kind of behavior occurs in the simulations of the LD function is that the phase space of the system is unbounded, and the action-based LD applied to a grid of initial conditions grows exponentially fast with the integration time due to the existence of a saddle equilibrium point at the origin. In this context, if one would like to reveal both manifolds in the same plot when the total LD (forward plus backward) is displayed, the approach that should be followed is to calculate what is known as the variable integration time LD \cite{GG2020a,ldbook2020}. The way that this variation of LDs is implemented is as follows. For any initial condition, the action-based LD is calculated along its trajectory for a given integration time $\tau$ or until the trajectory leaves a predefined region of the phase space, whatever happens first. We demonstrate in Fig, \ref{ld_saddle_varTime} that this procedure allows us to recover all the phase space structure in the same picture. To do so, we have integrated trajectories forward and backward in time for $\tau = 8$ or until they leave the square $[-8,8]\times[-8,8]$ centered about the saddle point at the origin, whatever happens first.

\section{The 1 DoF Harmonic Oscillator}
\label{sec:harmosc}

\begin{figure}[htbp]
	\begin{center}
    \includegraphics[scale=0.28]{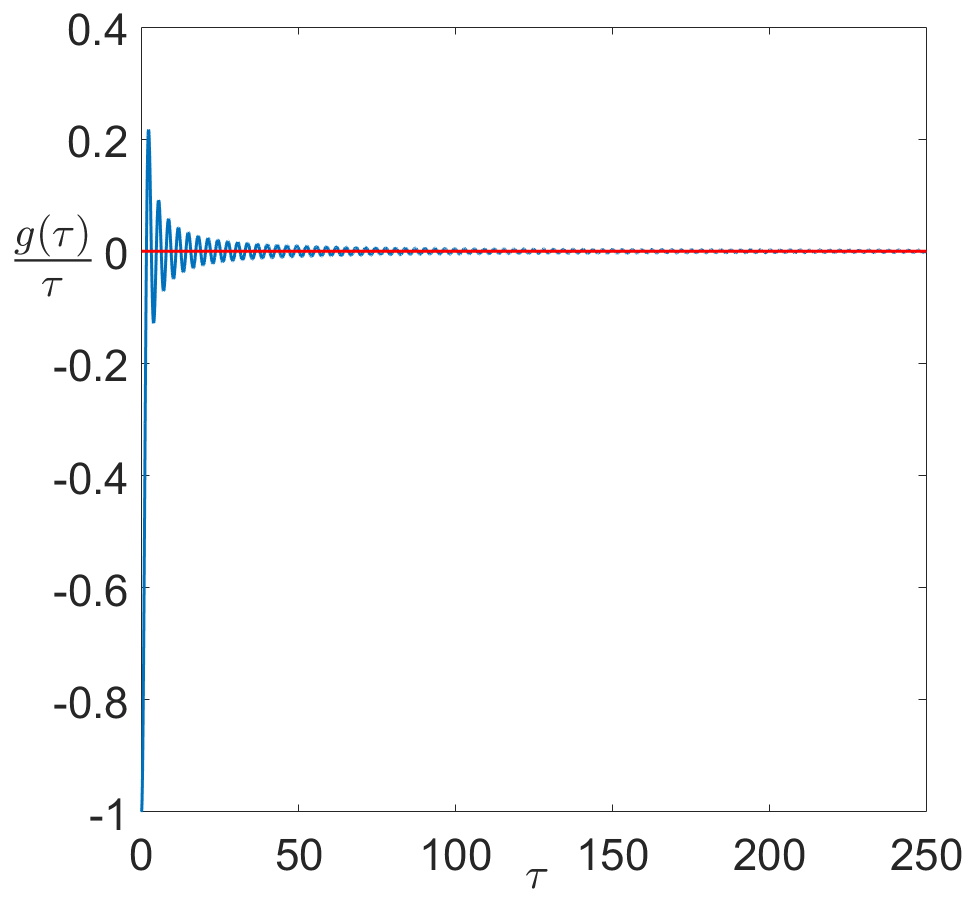}
	\end{center}
	\caption{Convergence of the long-time average of the forward LD, $\langle \mathcal{S}^{(f)}\rangle$, given in Eq. \eqref{func_conv}.}
	\label{conv_harmOsc}
\end{figure}

Here we will analytically show that the action-based LDs recovers the phase space of the simple harmonic oscillator. We will prove that the time-average of LDs converge to the tori that foliate the phase space of the system. Consider the 1 DoF Hamiltonian:
\begin{equation}
	\mathcal{H}(q,p) = \dfrac{p^2}{2m} + \dfrac{1}{2}m\omega^2q^2 \;,
\end{equation}
whose dynamics are governed by Hamilton's equations of motion:
\begin{equation}
	\begin{cases}
		\dot{q} = \dfrac{\partial \mathcal{H}}{\partial p} = \dfrac{p}{m} \\[.4cm]
		\dot{p} = -\dfrac{\partial \mathcal{H}}{\partial q} = -m \omega^2 q 
	\end{cases} \;.
\end{equation}
The solution to the harmonic oscillator are given by:
\begin{equation}
    q(t) = A \cos(\omega t) \quad,\quad p(t) = -m\omega A \sin(\omega t) \quad,\quad A = \dfrac{1}{\omega}\sqrt{\dfrac{2\mathcal{H}_0}{m}} \;,
\end{equation}
where $\mathcal{H}_0$ is the energy of the initial condition $\mathbf{x}_0 = (q_0,p_0) = (A,0)$. We calculate next the action-based LD in forward time:
\begin{equation}
    \mathcal{S}^{(f)}(\mathbf{x}_0,\tau) = \int_{0}^{\tau} \dfrac{p^2}{m} \, dt = m\omega^2A^2 \int_{0}^{\tau} \sin^2(\omega t) \, dt = \mathcal{H}_0 \left(\tau - \dfrac{\sin(2\omega \tau)}{2\omega}\right)
\end{equation}
Therefore, if we compute the time average of LDs:
\begin{equation}
\langle \mathcal{S}^{(f)}\rangle (\mathbf{x}_0,\tau)  = \dfrac{\mathcal{S}^{(f)}(\mathbf{x}_0,\tau)}{\tau}
\end{equation}
and take the limit as $\tau$ gets large we obtain:
\begin{equation}
 \mathcal{S}_{\infty} = \lim_{\tau \to \infty} \langle \mathcal{S}^{(f)}\rangle (\mathbf{x}_0,\tau) = \lim_{\tau \to \infty} \dfrac{\mathcal{S}^{(f)}(\mathbf{x}_0,\tau) }{\tau} = \lim_{\tau \to \infty} \mathcal{H}_0 \left(1 - \dfrac{\sin(2\omega \tau)}{2\omega\tau}\right) = \mathcal{H}_0   
\end{equation}
If we define:
\begin{equation}
    g(\tau) = \left(\langle \mathcal{S}^{(f)}\rangle -  \mathcal{S}_{\infty}  \right) \, \tau = -\dfrac{\mathcal{H}_0}{2\omega} \sin(2\omega\tau) 
    \label{func_conv}
\end{equation}
then the Fourier Transform of $g$ would yield twice the frequency of the harmonic oscillator, i.e. $2\omega$. Therefore, we can recover information of regular quasiperiodic motion of the system from the long-time averages of the LD function. In Fig. \ref{conv_harmOsc} we show how the time average of the forward component of LDs, that is, $\langle \mathcal{S}^{(f)}\rangle$ tends towards $ \mathcal{S}_{\infty}$ as time $\tau$ goes to infinity.

\section{A Two DoF Hamiltonian Model}
\label{sec:2dof}

\begin{figure}[htbp]
	\begin{center}
		A)\includegraphics[scale=0.26]{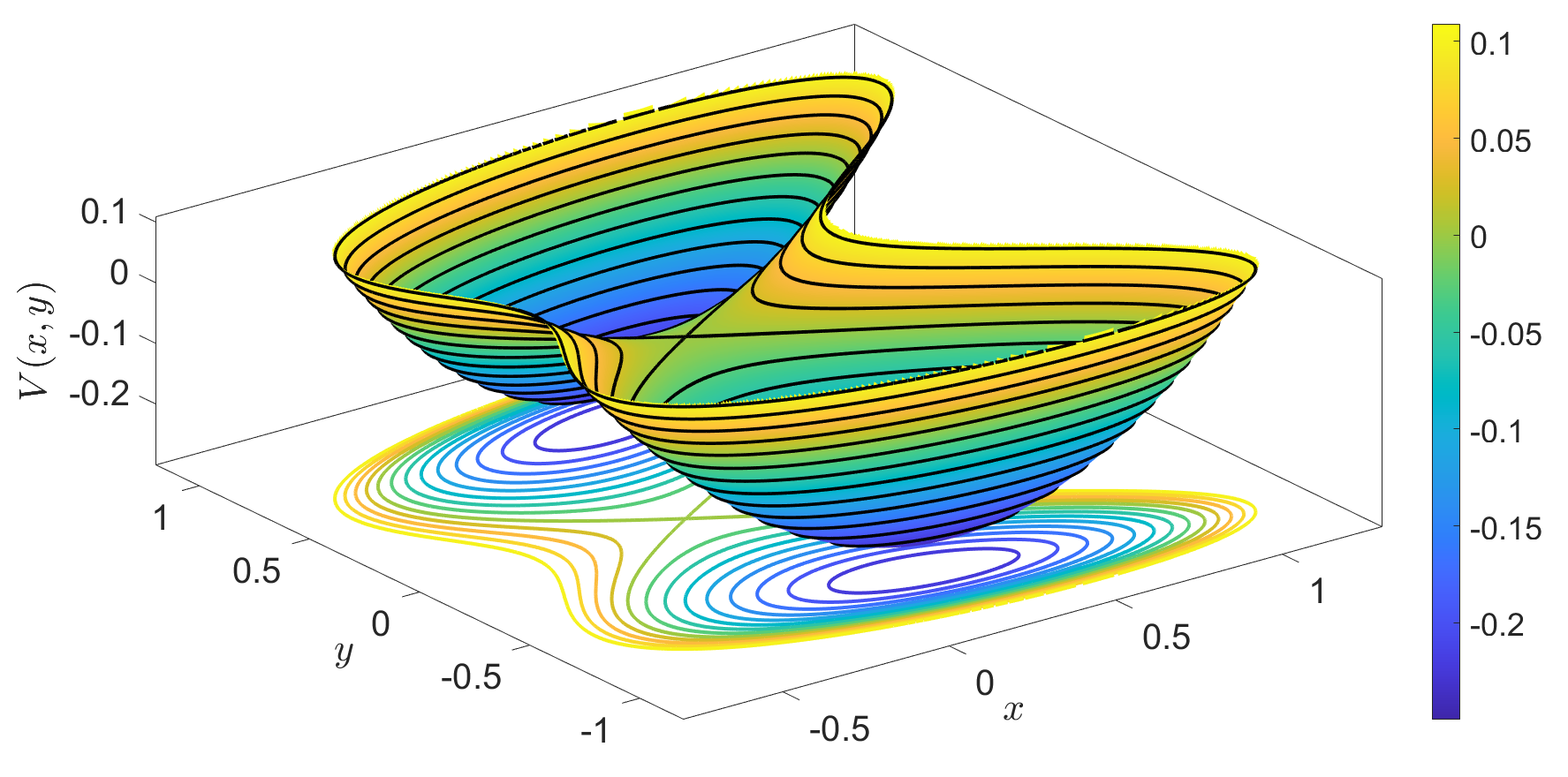}
		B)\includegraphics[scale=0.24]{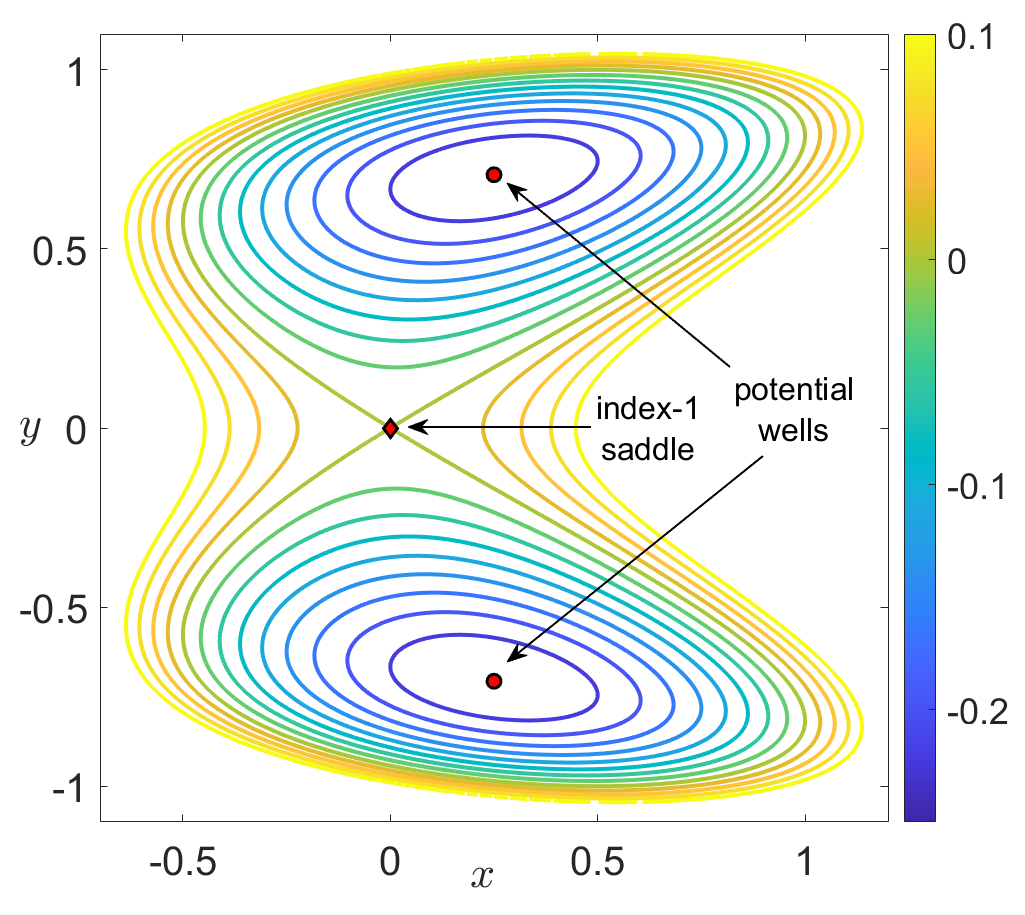}
		C)\includegraphics[scale=0.38]{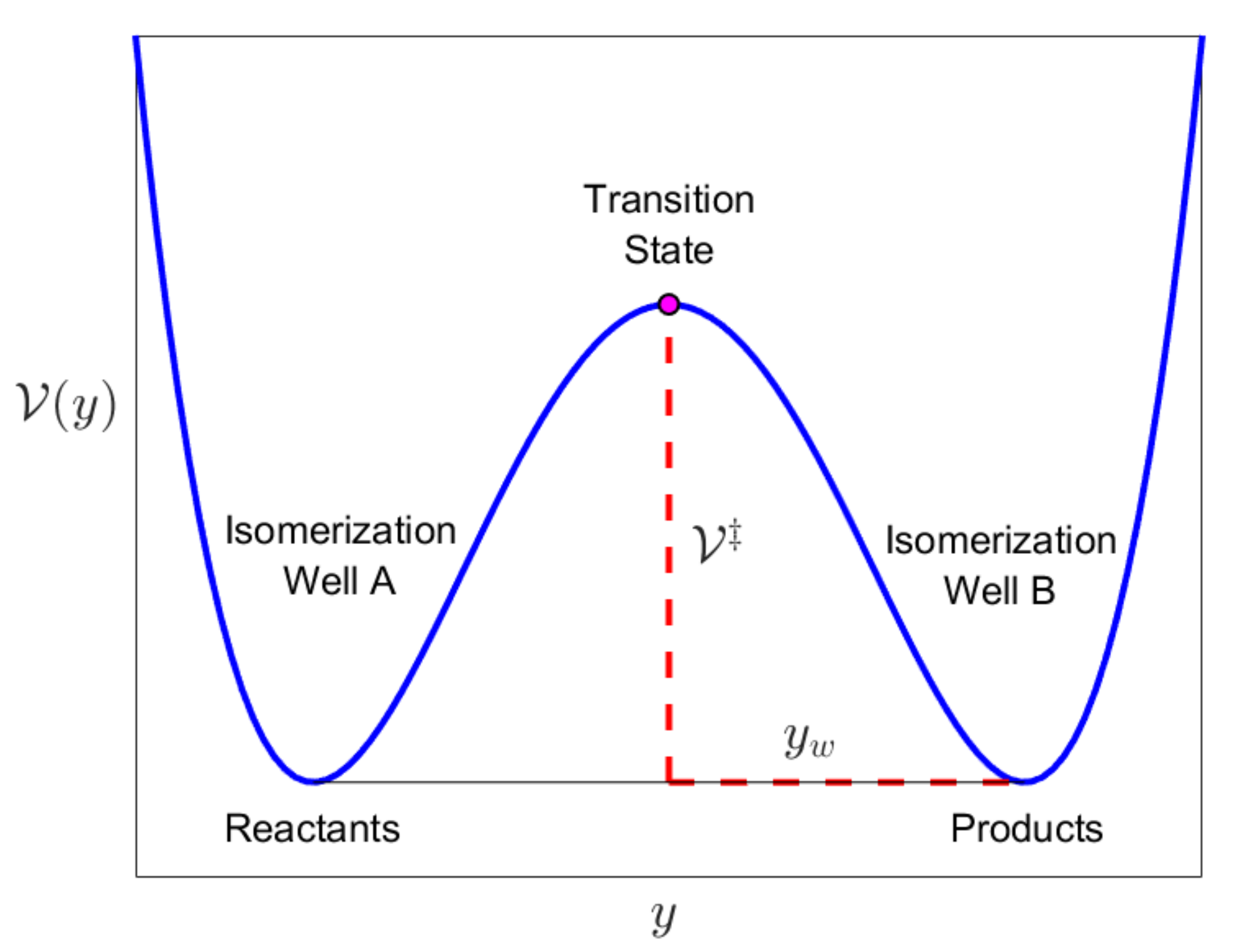}	
	\end{center}
	\caption{A) and B) Potential energy surface landscape and energy contours for the system parameters $\mathcal{V}^{\ddagger} = 1/4$, $y_w = \sqrt{2}/2$, $\omega = 1$ and $c = 1/2$. We have marked the location of the index-1 saddle at the origin with a red diamond, and the potential wells with red circles. C) Symmetric quartic bistable potential in the $y$ DoF showing the location of the isomerization wells and the transition state. The parameter $\mathcal{V}^{\ddagger}$ measures the potential barrier height measured from the bottom of any of the potential wells, and $y_w$ is the horizontal distance from the index-1 saddle to any of the wells.}
	\label{PES}
\end{figure}

In this section we show how the action-based LD can be used to reveal the phase space structures of a 2 DoF Hamiltonian system. In particular, we consider a 2 DoF Hamiltonian model that has been used to study proton transfer and isomerization reactions in chemistry \cite{paz95}. The Hamiltonian function is given as the sum of kinetic plus potential energies of the form:
\begin{equation}
	\mathcal{H}(x,y,p_x,p_y) = T(p_x,p_y) + V(x,y) = \dfrac{p_x^2}{2m} + \dfrac{p_y^2}{2m} + \dfrac{\mathcal{V}^{\ddagger}}{y_w^4} \, y^2 \left(y^2 - 2y_w^2\right) + \dfrac{1}{2} m \omega^2 \left(x - \dfrac{c y^2}{m\omega^2}\right)^2
\end{equation}
where $m$ is the mass of the $x$ and $y$ DoFs. The potential energy surface of the system consists of a symmetric double well potential in the $y$ DoF, with a saddle at the origin, and two wells located at $\pm y_w$. The barrier height is given by the value of $\mathcal{V}^{\ddagger} > 0$. The $y$ DoF is coupled to the $x$ DoF, which represents a bath coordinate (or harmonic oscillator), by means of a quadratic coupling, and $\omega$ is the angular frequency of the $x$ DoF. The strength of the coupling is defined by the constant $c$. In Fig. \ref{PES} we depict the topography and energy contours of the PES, together with schematic representation of the symmetric double well potential in the $y$ DoF.

Hamilton's equations of motion for this system are given by:
\begin{equation}
	\begin{cases}
		\dot{x} = \dfrac{\partial \mathcal{H}}{\partial p_x} = \dfrac{p_x}{m} \\[.5cm]
		\dot{y} = \dfrac{\partial \mathcal{H}}{\partial p_y} = \dfrac{p_y}{m} \\[.5cm]
		\dot{p}_x = -\dfrac{\partial \mathcal{H}}{\partial x} = -m\omega^2 x + cy^2 \\[.4cm]
		\dot{p}_y = -\dfrac{\partial \mathcal{H}}{\partial y} = 2y \left[\left(\dfrac{2\mathcal{V}^{\ddagger}}{y_w^2} + cx\right) - \left(\dfrac{2\mathcal{V}^{\ddagger}}{y_w^4} + \dfrac{c^2}{m\omega^2}\right)y^2 \right]
	\end{cases} \;,
	\label{hameq}
\end{equation}
and for simplicity in our analysis we will set $m = 1$. The phase space of this dynamical system is four-dimensional, and since energy is conserved, motion is constrained to a three-dimensional energy hypersurface. It is straightforward to show that this system has three equilibrium points, one is an index-1 saddle at the origin with eigenvalues:
\begin{equation}
	\lambda_{1,2} = \pm \dfrac{2}{y_w} \sqrt{\dfrac{\mathcal{V}^{\ddagger}}{m}} \quad,\quad \lambda_{3,4} = \pm \omega \, i
\end{equation}
and two centers symmetrically located with respect to the $x$-axis at the points $\left(cy_w^2/(m\omega^2),\pm y_w,0,0\right)$, which indicate the bottom of the potential wells. For energies above the energy of the index-1 saddle it is possible for trajectories to cross between the wells. The phase space structures responsible for mediating this transport are the stable and unstable manifolds associated with an unstable periodic orbit, which is guaranteed to exist for a range of energies above the energy of the index-1 saddle by the Lyapunov subcenter theorem \cite{moser1958,lyapunov1992general}.

We now show that the action-based LDs can  reveal these phase space structures that characterize transport for the Hamiltonian system in Eq. \eqref{hameq}.  We begin by fixing the energy of the system to $\mathcal{H} = 0.1$ and set the values for the model parameters $\mathcal{V}^{\ddagger} = 1/4$, $y_w = \sqrt{2}/2$, $\omega = 1$ and $c = 1/2$. We will use the action-based LDs to show how these structures appear in  the surfaces of section:
\begin{equation}
\Sigma_1 = \left\lbrace (x,y,p_x,p_y) \in \mathbb{R}^4 \; \Big| \; x = 0 \;,\; p_x \geq 0 \right\rbrace \quad,\quad \Sigma_2 = \left\lbrace (x,y,p_x,p_y) \in \mathbb{R}^4 \; \Big| \; y = -y_w \;,\; p_y \geq 0 \right\rbrace
\label{sos}
\end{equation}
where $\Sigma_1$ contains the index-1 saddle at the origin, and $\Sigma_2$ is located at the bottom of the well of the potential energy surface.

\begin{figure}[htbp]
	\begin{center}
	A)\includegraphics[scale=0.28]{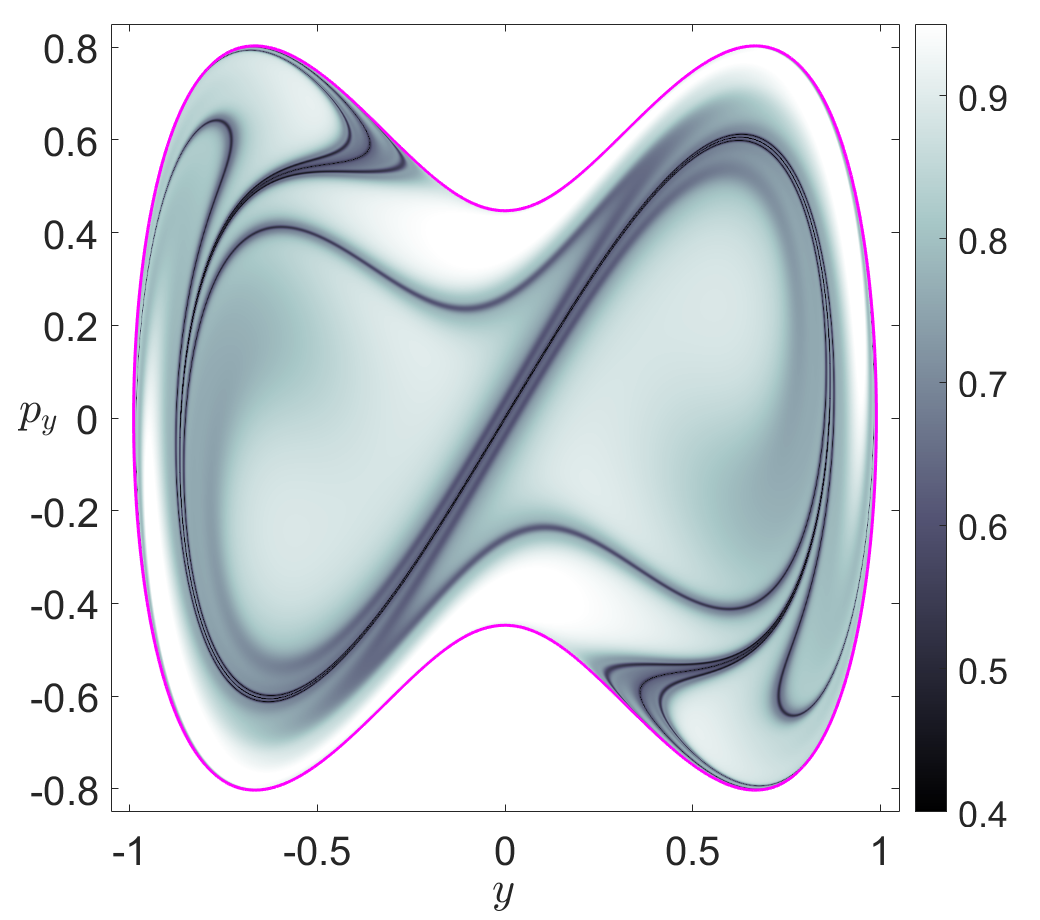}
	B)\includegraphics[scale=0.28]{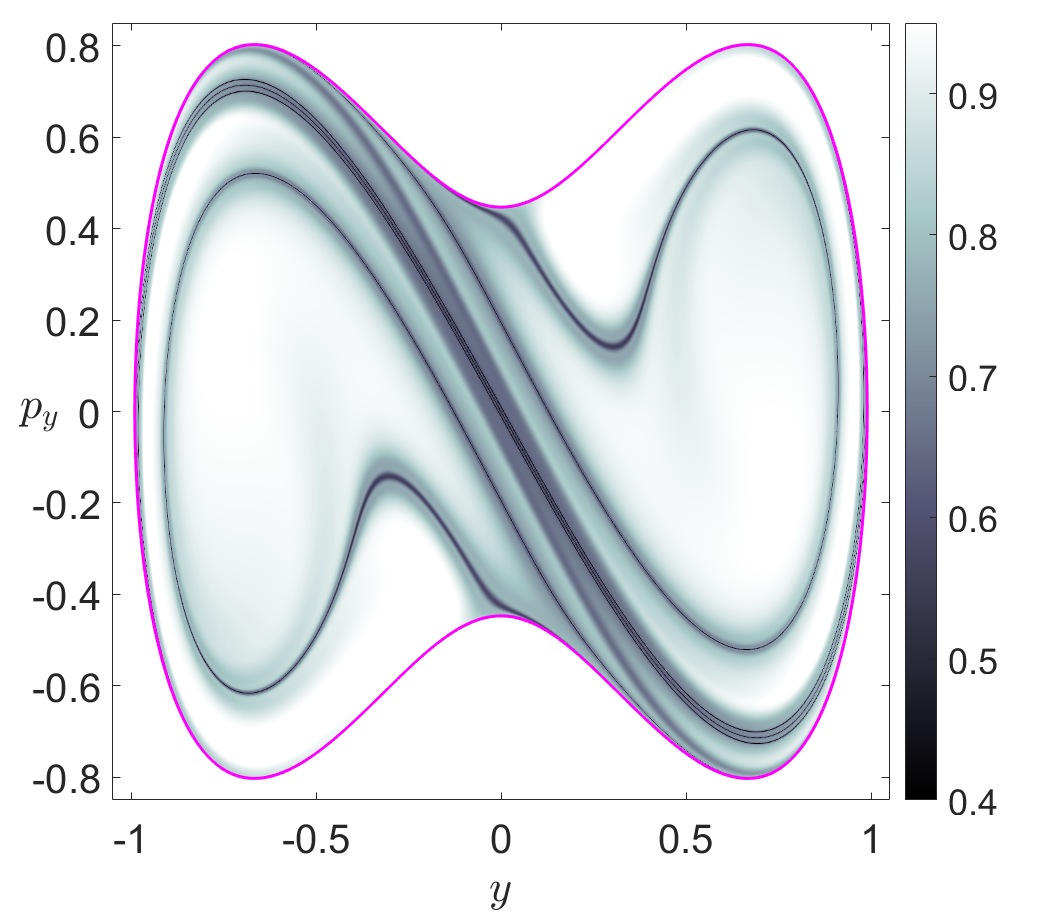}
	C)\includegraphics[scale=0.28]{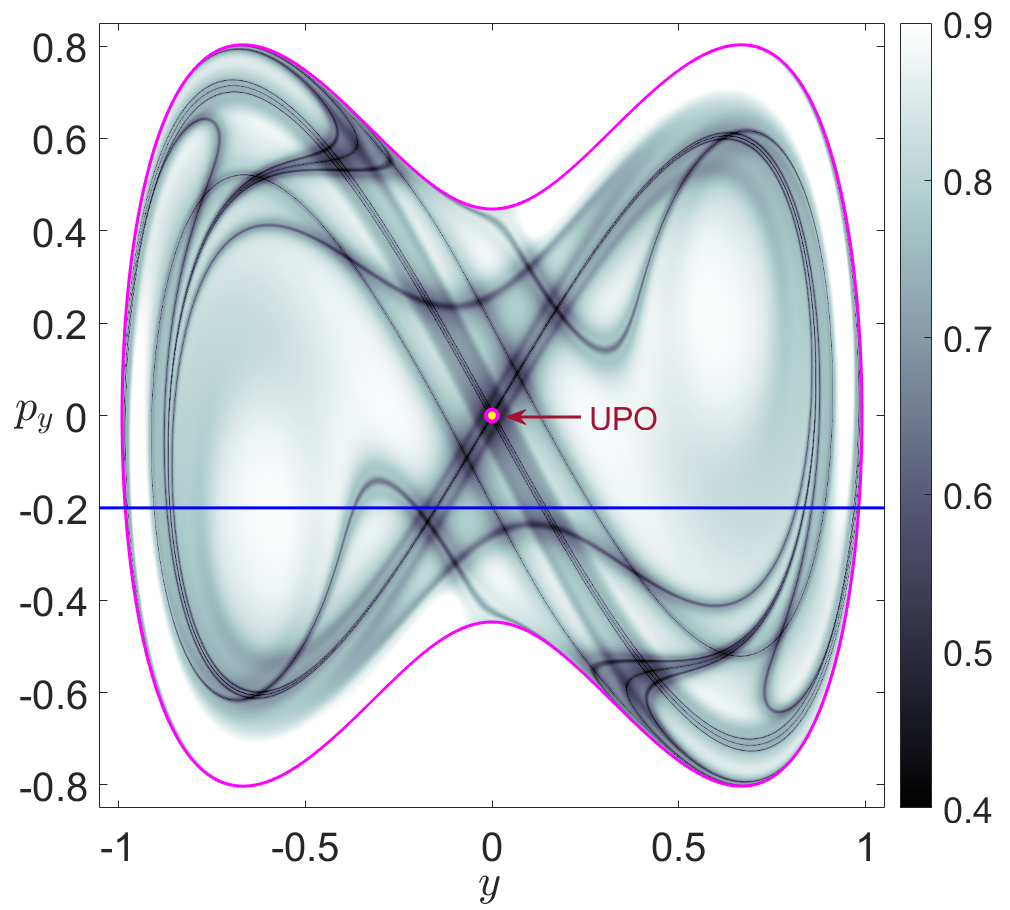}
	D)\includegraphics[scale=0.28]{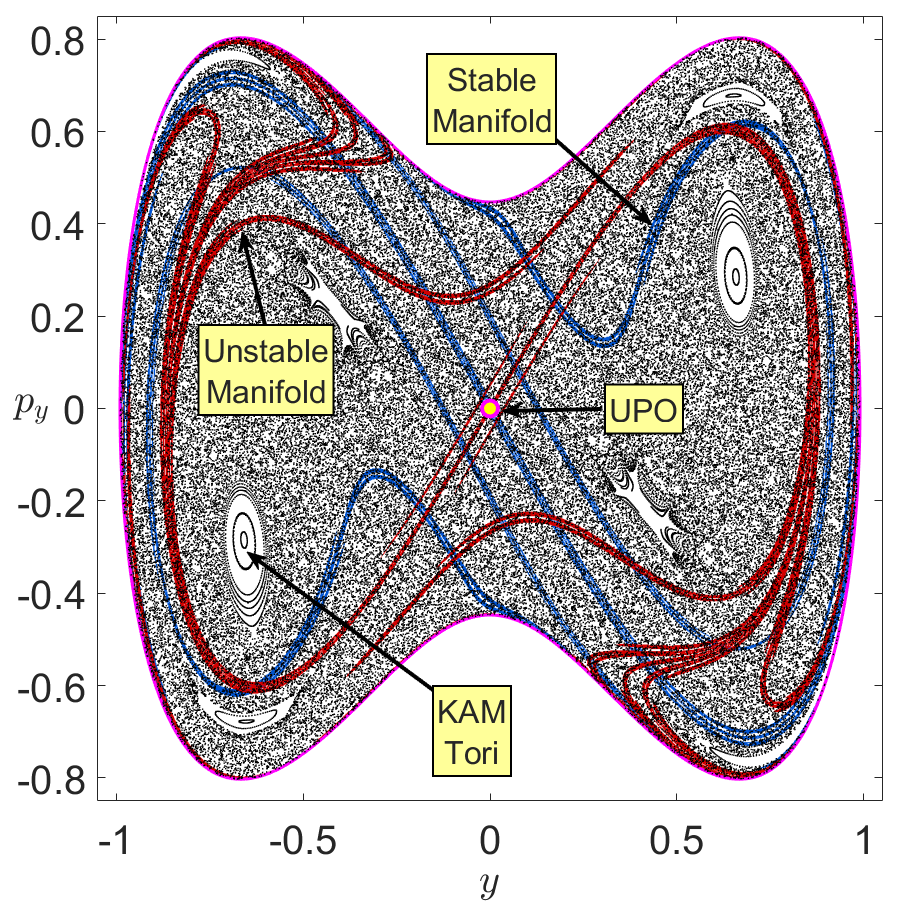}
	E)\includegraphics[scale=0.32]{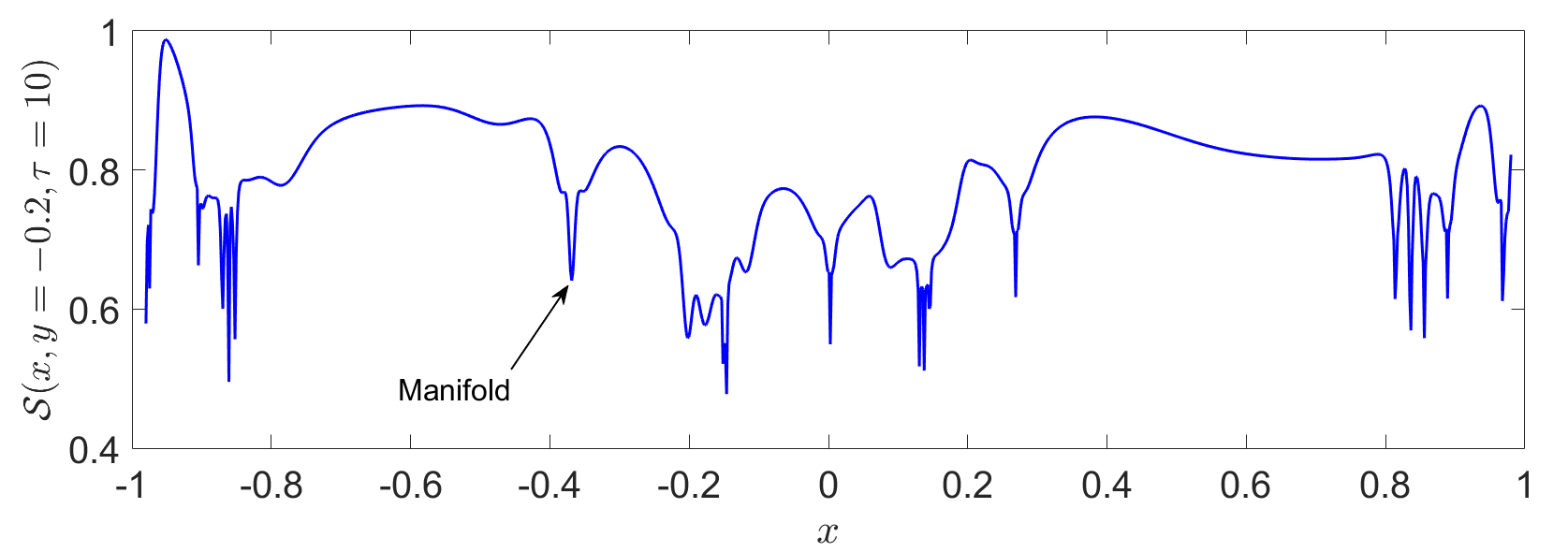}	
	\end{center}
	\caption{Phase space of the Hamiltonian system at energy $\mathcal{H} = 0.1$, as revealed by Lagrangian descriptors with $\tau = 10$ on the section $\Sigma_1$ in Eq. \eqref{sos}. A) Backward component of LDs; B) Forward component of LDs; C) Total LD (sum of forward plus backward components); D) Poincar\'e section superimposed with the stable (blue) and unstable (red) manifolds extracted from the LD scalar field presented in C); E) LD values along the line $p_y = -0.2$ depicted in panel C) to show how the LD scalar field displays 'singular features' at the manifold locations.}
	\label{psec_x0}
\end{figure}

In Fig. \ref{psec_x0} we show the action-based LDs (computed for $\tau =10$) and Poincar\'e maps on $\Sigma_1$. In panel A we see that the action-based backward LD highlights the unstable manifold of the UPO (which is indicated in the figure as a circle) and in panel B we see that the action-based forward LD highlights the stable manifold of the UPO. Panel C shows the forward and backward LDs together. In panel D the stable and unstable manifolds are extracted from the LD plots and shown with a Poincar\'e map. Panel E shows how the 'singularities' of the LD capture the manifolds. It is important to remark here that the LD scalar field attains a local minimum at the location of the UPO, and this important property of the action-based integral formulation of LDs can be conveniently used to detect these objects and implement a minimization algorithm to compute them. These periodic trajectories play an essential role for the study of transition phenomena across the phase space bottleneck that exists in the neighborhood of index-1 saddle between both wells, since they provide the scaffolding for constructing a dividing surface with the local non-recrossing and minimal flux properties that 'reactive' trajectories cross in their path from one well to the other \cite{waalkens2004direct}. These phase space structures are fundamental for the development of Transition State Theory in chemistry \cite{wiggins2001,uzer2002geometry} or in the design of space missions in astrodynamics \cite{Koon2011}.

\begin{figure}[htbp]
	\begin{center}
		A)\includegraphics[scale=0.28]{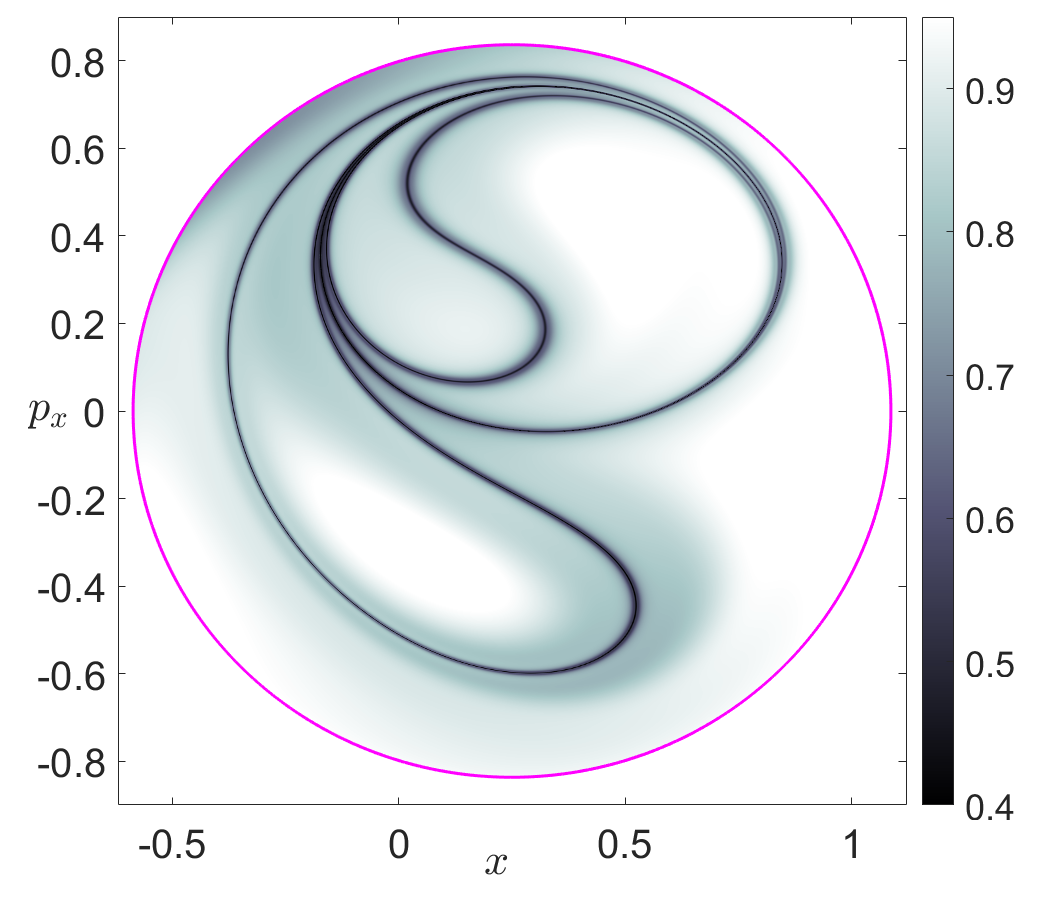}
		B)\includegraphics[scale=0.28]{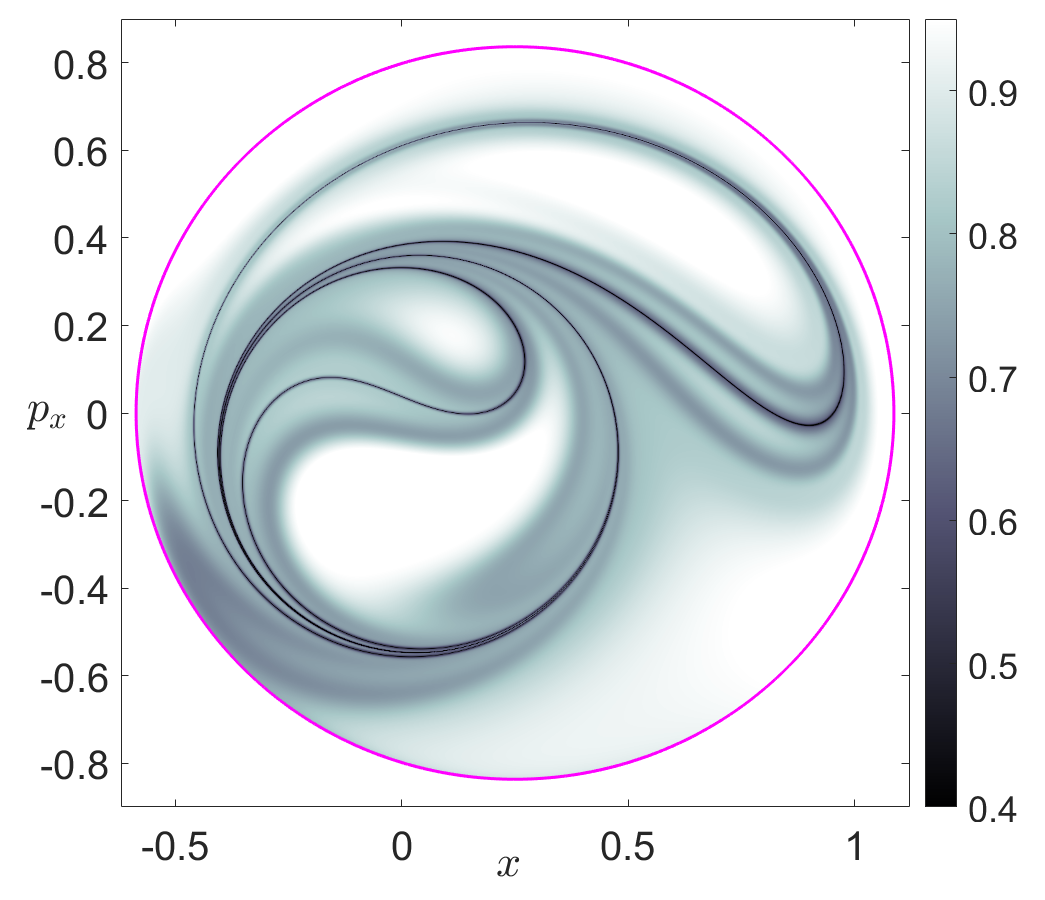}
		C)\includegraphics[scale=0.28]{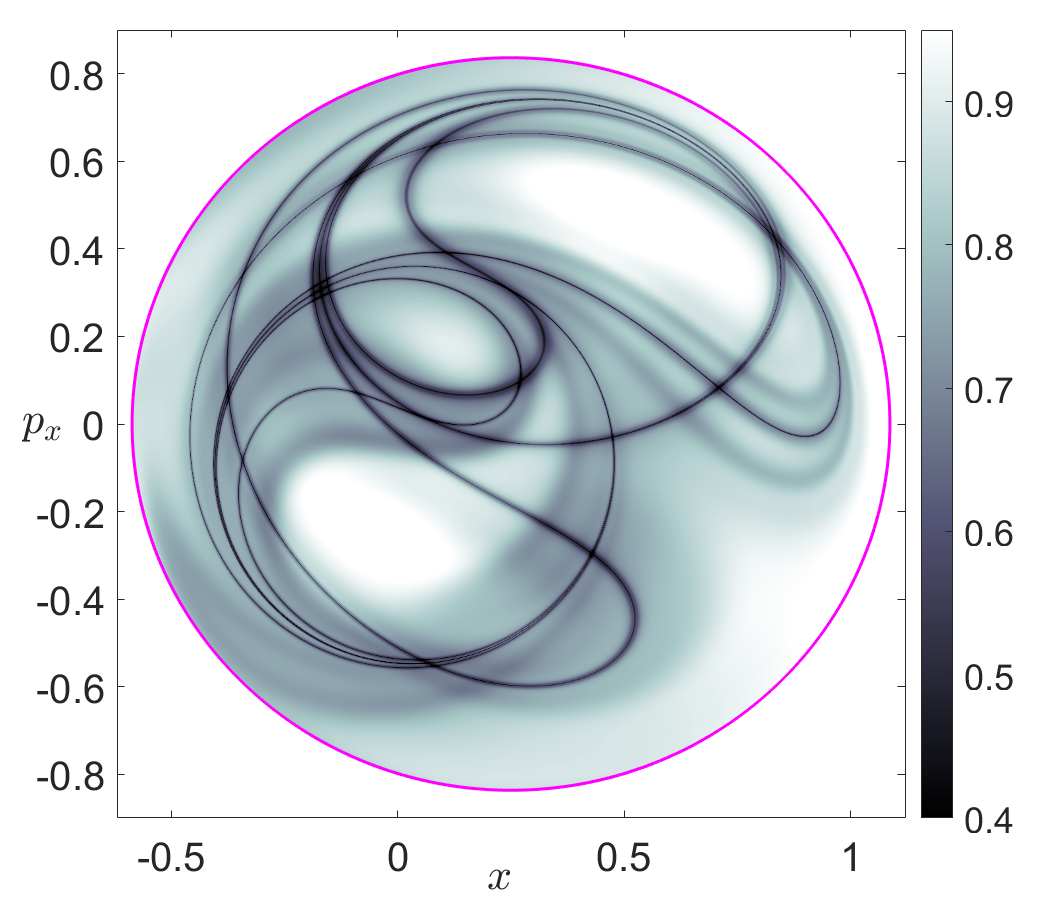}
		D)\includegraphics[scale=0.28]{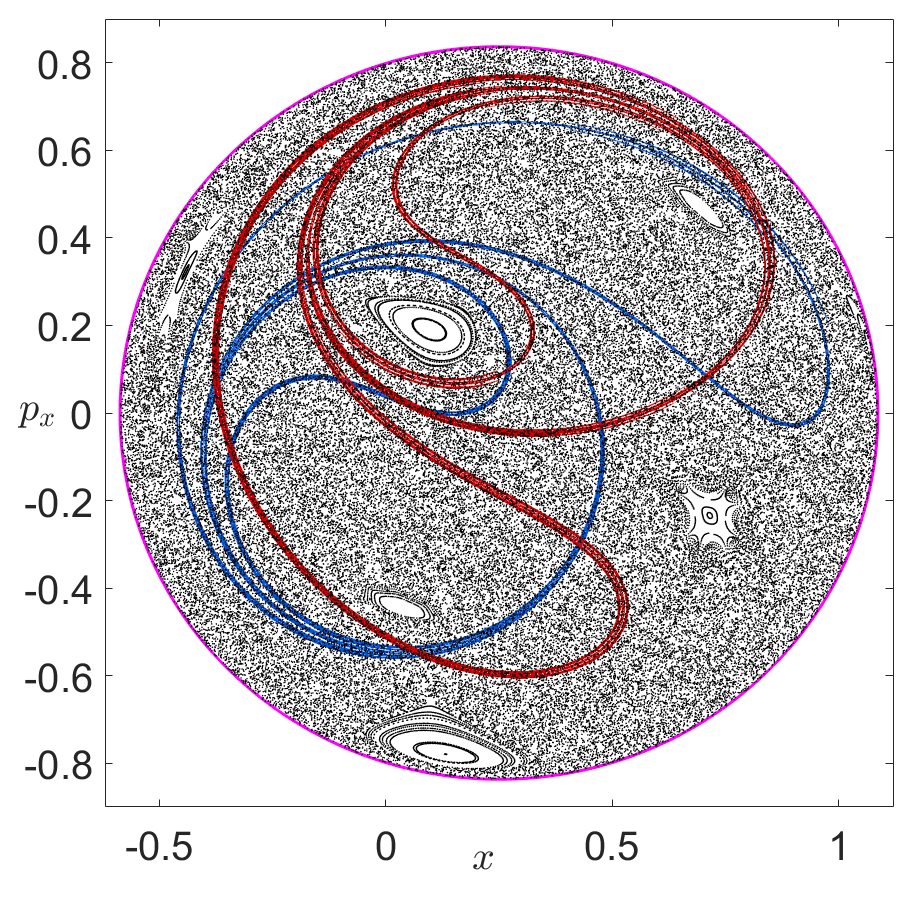} 	
	\end{center}
	\caption{Phase space of the Hamiltonian system at energy $\mathcal{H} = 0.1$, as revealed by Lagrangian descriptors with $\tau = 10$ on the section $\Sigma_2$ in Eq. \eqref{sos}. A) Backward component of LDs; B) Forward component of LDs; C) Total LD (sum of forward plus backward components); D) Poinca\'e section superimposed with the stable (blue) and unstable (red) manifolds extracted from the LD scalar field shown in C).}
	\label{psec_ywell}
\end{figure}

We display in Fig. \ref{psec_ywell} the action-based LDs and Poincar\'e maps on $\Sigma_2$. Similar panels as those displayed in Fig. \ref{psec_x0} are used, but in this case the section $\Sigma_2$ shows the extent to which the manifolds can penetrate into the well and how they give rise to transport into and out of the phase space region corresponding to this well.

\begin{figure}[htbp]
	\begin{center}
	A)\includegraphics[scale=0.29]{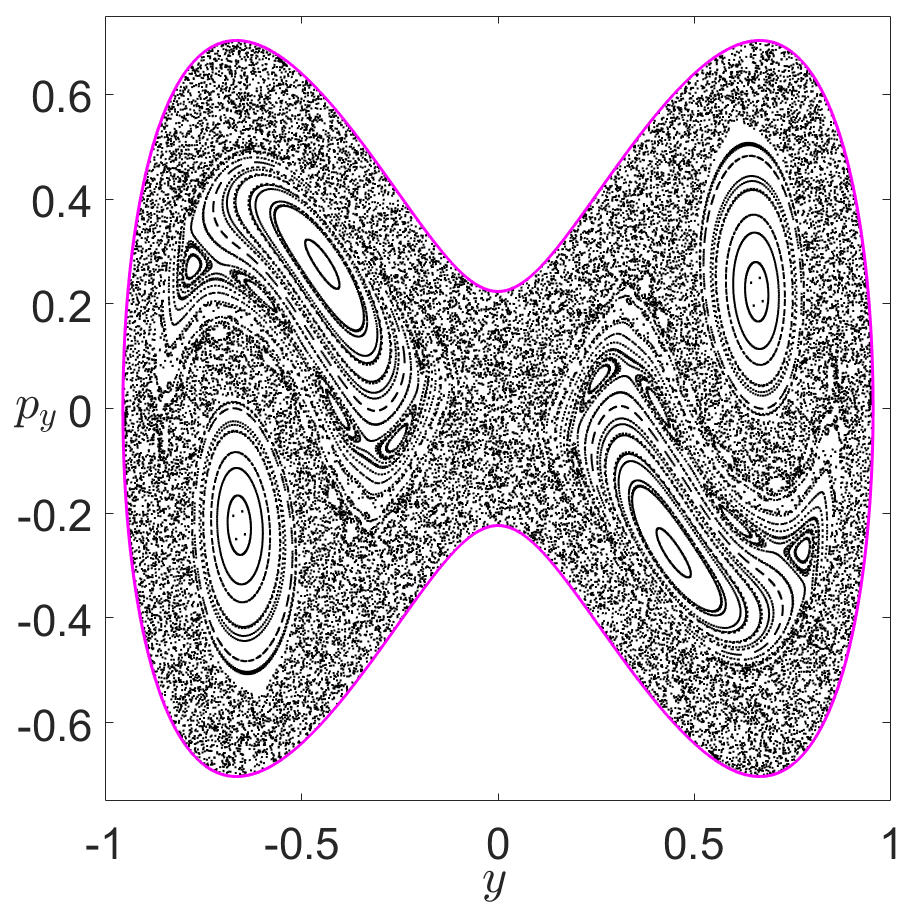}
	B)\includegraphics[scale=0.29]{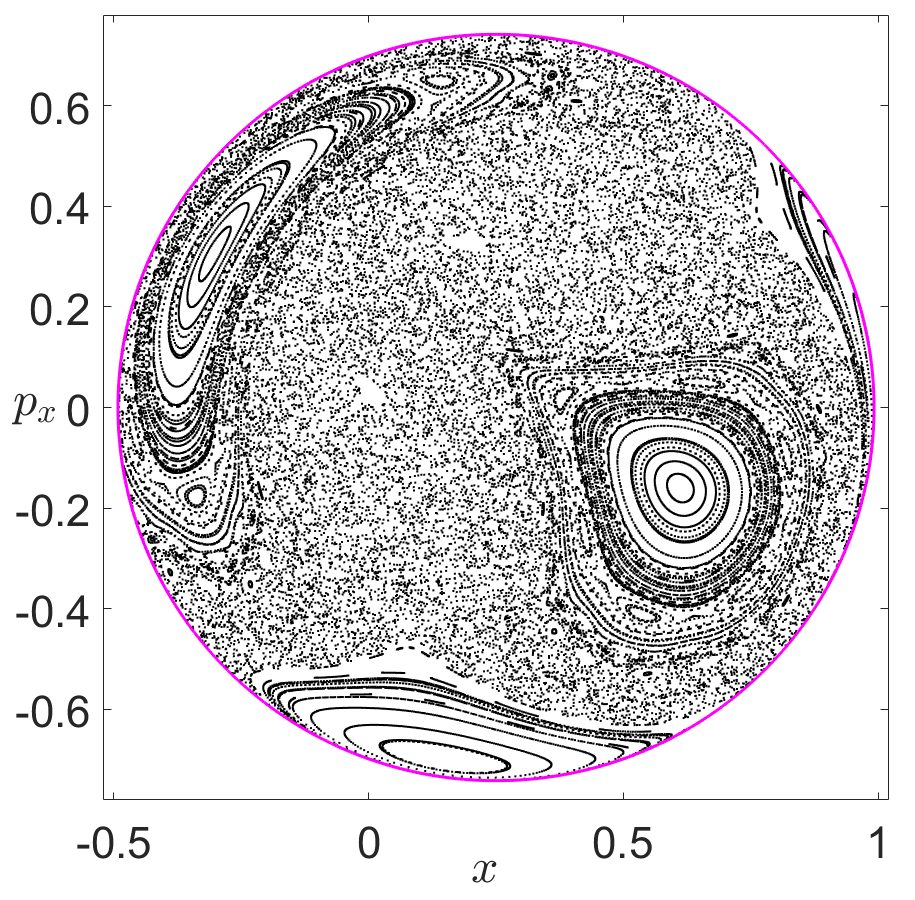}
	C)\includegraphics[scale=0.26]{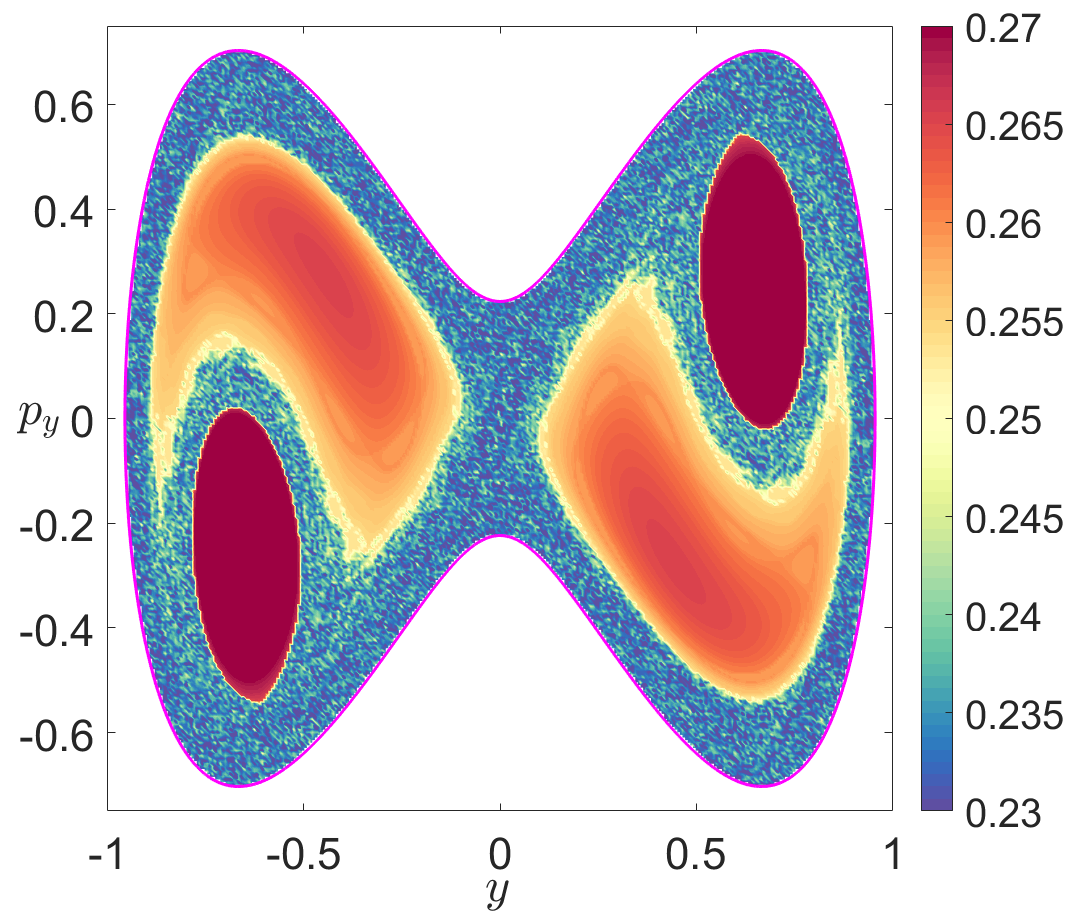}
	D)\includegraphics[scale=0.26]{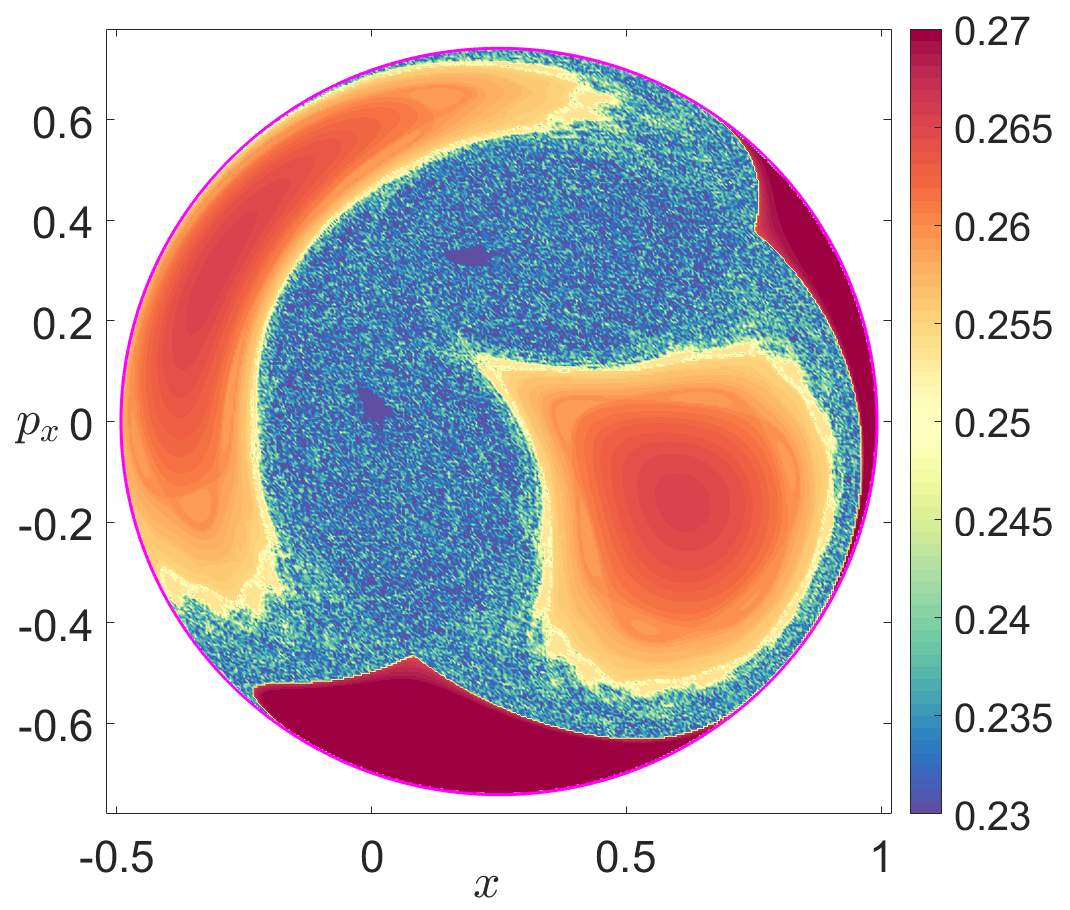}
	E)\includegraphics[scale=0.26]{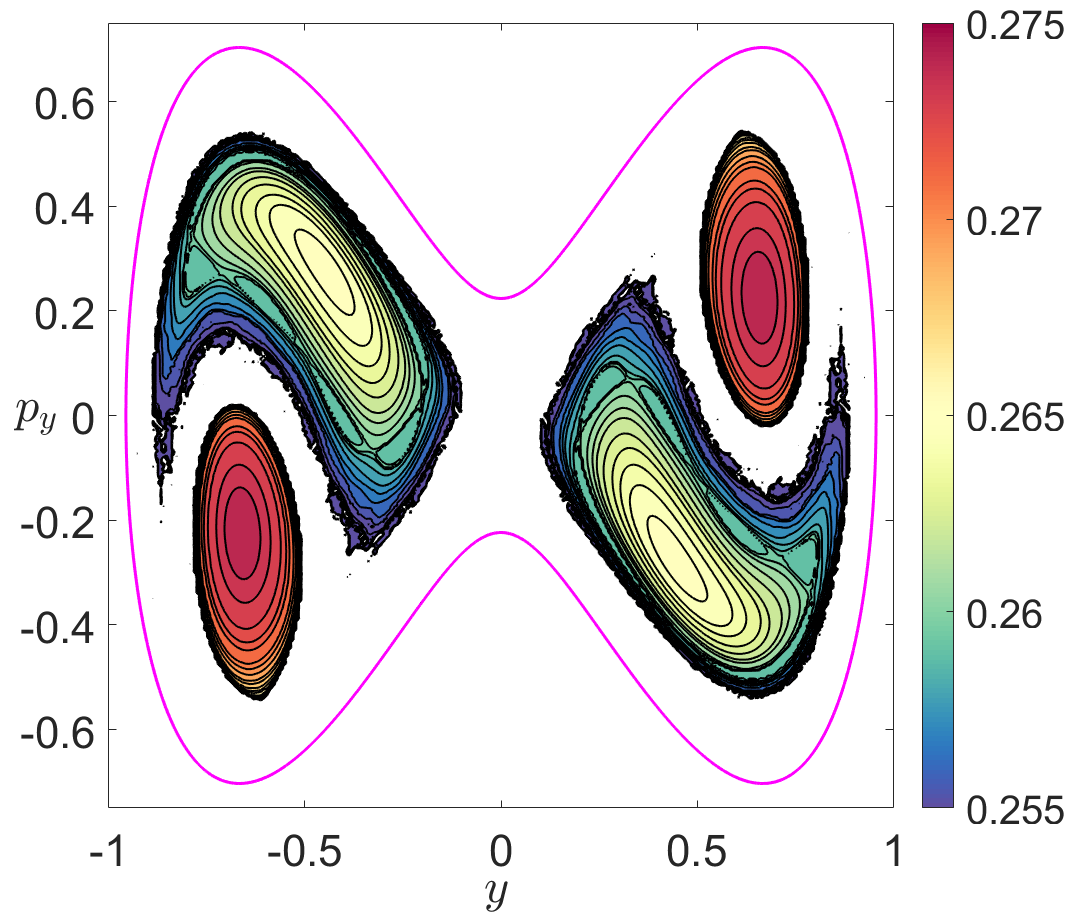}
	F)\includegraphics[scale=0.26]{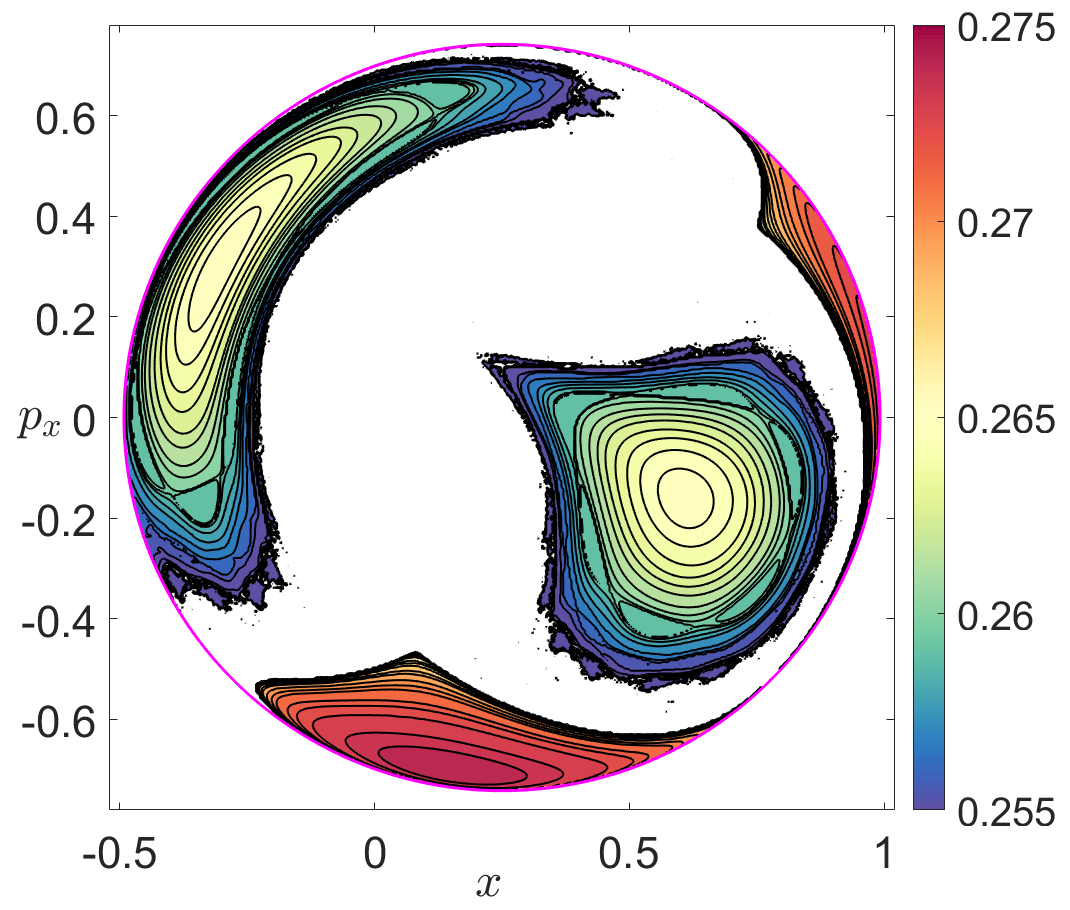}
	\end{center}
	\caption{Poincar\'e sections and time average of Lagrangian descriptors calculated on the sections $\Sigma_1$ (left column) and $\Sigma_2$ (right column) for the Hamiltonian system in Eq. \eqref{hameq} with energy $\mathcal{H} = 0.025$. A) and B) Poincar\'e map in forward time for $t = 1500$. C) and D) Time average of the total LD (forward plus backward) for an integration time $\tau = 750$. E) and F) Filled contours for the time averaged LD scalar field.}
	\label{tavg}
\end{figure}

We finish this section by demonstrating how LDs defined in terms of the action integral can be used to recover the KAM tori characterizing regular quasiperiodic motion for the Hamiltonian system in Eq. \eqref{hameq}. This type of analysis, based on the computation of the long-term time averages of the LD scalar field was first introduced in \cite{lopesino2017} for simple two and three dimensional systems using the $p$-norm definition of LDs, and developed further in the  three-dimensional setting in \cite{gg2018}. This property of the LD method was rigorously proven for the case $p = 1$ using the harmonic oscillator in \cite{lopesino2017}, and later it was extended to other values of the $p$-norm in \cite{naik2019a}. This approach has also been addressed for maps in \cite{GG2020}. Recently, it has been reproduced in \cite{montes2021}, together with the original derivation of the case $p = 1$ already discussed in \cite{lopesino2017}, to study regular motion in the H\'enon-Heiles potential. In this paper, we apply the action-based definition of LDs in Eq. \eqref{eq:LD_S} to analyze the phase space of the 2 DoF Hamiltonian system in Eq. \eqref{hameq} with an energy $\mathcal{H}_0 = 0.025$, which is above that of the index-1 saddle at the origin. As we proved in Sec. \ref{sec:harmosc}, long-term time averages of LDs converge to invariant tori, and these structures can be extracted directly from the LD scalar field. To illustrate this capability we compute LDs on the sections $\Sigma_1$ and $\Sigma_2$ given in Eq. \eqref{sos} for an integration time $\tau = 750$ both forward and backward. The LD output is averaged in time as follows:
\begin{equation}
\langle \mathcal{S}\rangle (\mathbf{x}_0,\tau)  = \dfrac{\mathcal{S}^{(f)}(\mathbf{x}_0,\tau) + \mathcal{S}^{(b)}(\mathbf{x}_0,\tau)}{\tau}
\end{equation}
and the results of these simulations are displayed in Fig. \ref{tavg}. As a validation, we compare the KAM tori obtained from the contours of the time-averaged LD values with those visualized directly from a classical Poincar\'e map calculated by running an ensemble of initial conditions forward in time for $t = 1500$ units.

\section{The Stochastically Forced Duffing Oscillator}
\label{sec:stochsys}

In this section we illustrate how the action-based LD can be used to analyze the phase space of a random dynamical system \cite{duan2015}. We do so by looking at the dynamics of the stochastically forced Duffing oscillator, which gives a simple model to explore the rate at which a Brownian particle escapes from a potential well over a potential barrier. This is known in the literature as the Kramers' problem \cite{melnikov91} and has many applications in in chemical physics and also in biophysics. A similar type of analysis was carried out with LDs in \cite{balibrea2016lagrangian} using an alternative definition based on the $p$-norm of the vector field that determines the dynamical system under study.

\begin{figure}[htbp]
	\begin{center}
	A)\includegraphics[scale=0.28]{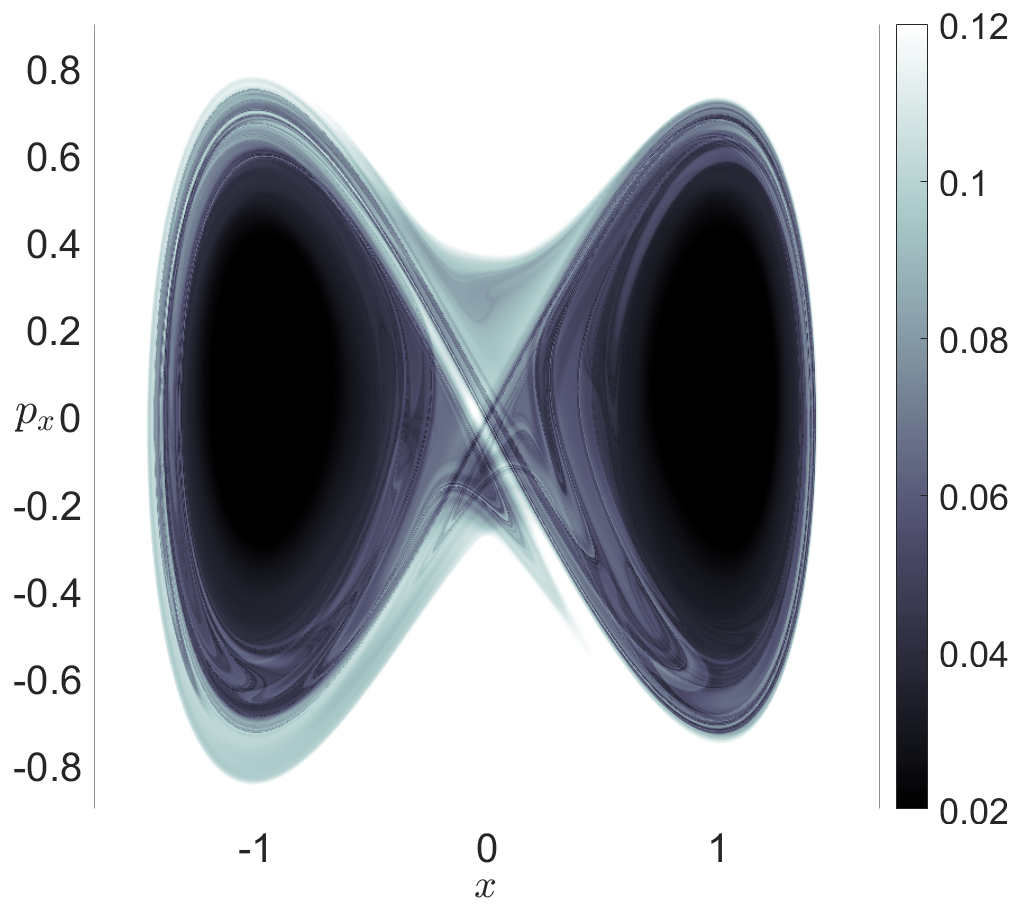}
	B)\includegraphics[scale=0.28]{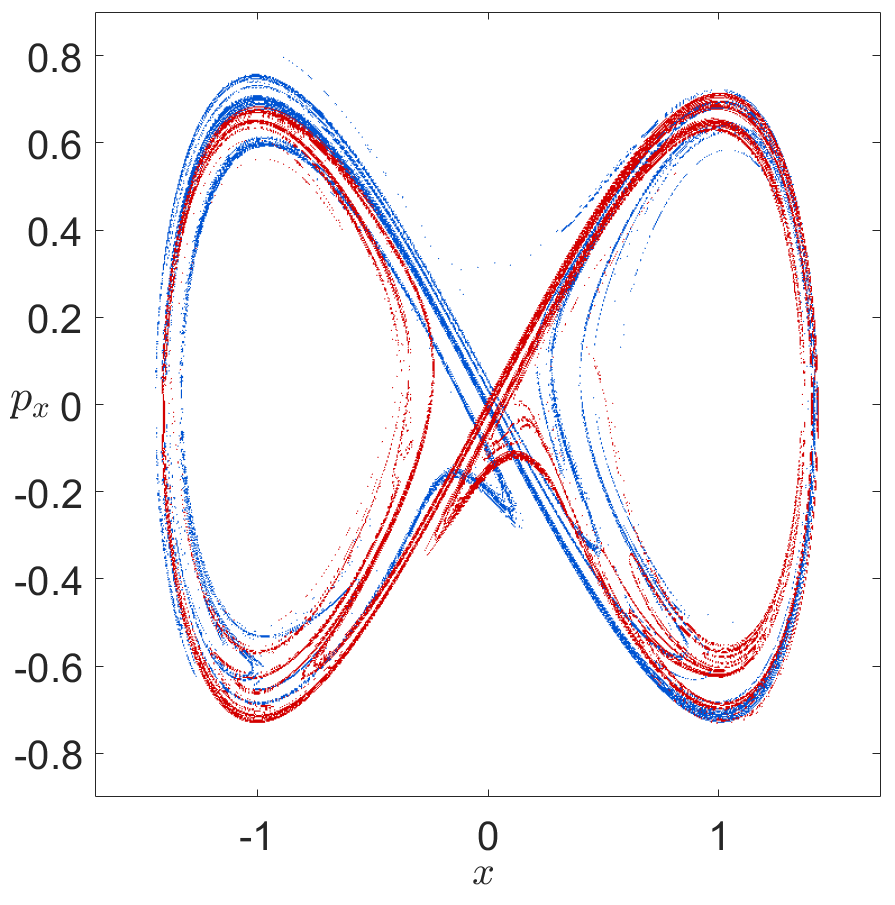}
	C)\includegraphics[scale=0.28]{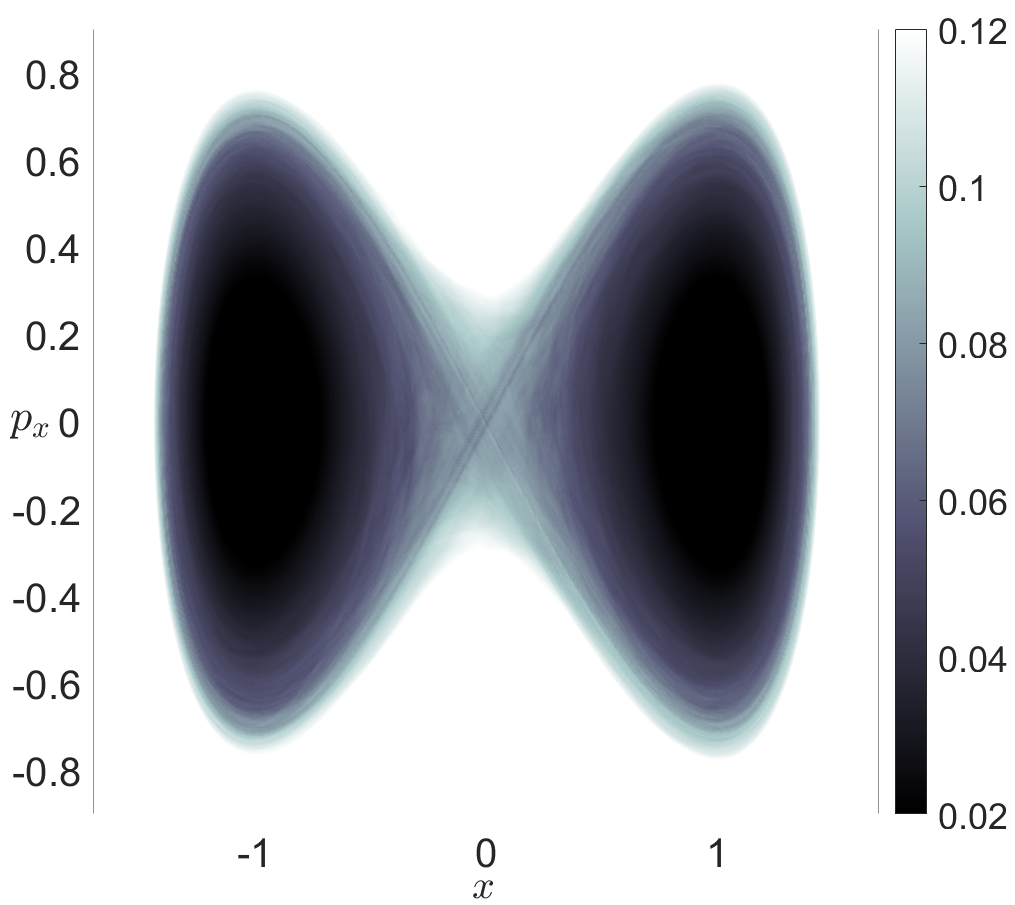}
	D)\includegraphics[scale=0.28]{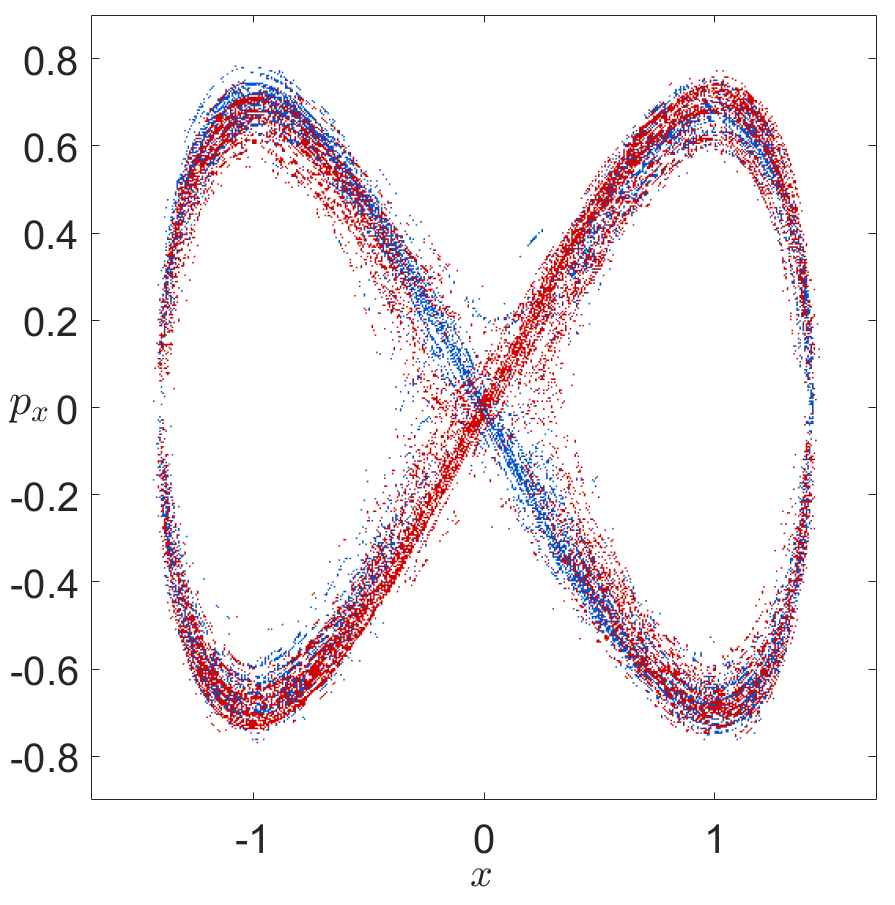} 	
	\end{center}
	\caption{Phase space of the stochastic Duffing system in Eq. \eqref{stoch_duff}, as revealed by the action-based Lagrangian descriptor with $\tau = 35$. A) $N = 1$ random experiment; C) Average LD for $N = 25$ random experiments. C) and D) correspond to the stable (blue) and unstable (red) manifolds extracted from the LD scalar field shown in panels A) and B), respectively.}
	\label{phase_duff}
\end{figure}

The stochastic Duffing system is described by the random dynamical system:
\begin{equation}
	\begin{cases}
		d X_t = Y_t \, dt \\[.1cm]
		d Y_t = \left(X_t - X_t^3\right) dt + \sigma \, dW_t 
	\end{cases} \,,
	\label{stoch_duff}
\end{equation}
where $W_t$ is a two-sided Wiener process and $\sigma$ represents the stochastic forcing strength which we consider to be additive, that is, it does not depend on the state of the system. In particular, we will give it a constant value $\sigma = 0.025$. To study the phase portrait of this SDE, we apply Lagrangian descriptors at an initial time $t = 0$ to a grid of $600\times600$ initial conditions on the rectangle $[-1.7,1.7]\times[-0.9,0.9]$, integrating trajectories forward and backward for a time $\tau = 35$. To do so, we have implemented the Euler-Maruyama scheme \cite{kloeden1992stochastic} using a time step of $\Delta t = 0.005$. First, we carry out only one experiment and the results of this numerical simulation are shown in panels A) and B) of Fig. \ref{phase_duff}. Our next goal is to reconstruct the most likely phase portrait for the Duffing oscillator by means of performing $25$ simulations, and taking the average of the resulting LD scalar fields. The outcome of this calculation is displayed in Fig. \ref{phase_duff} C) and D). Notice how the method nicely highlights the homoclinic tangle formed by the stable and unstable manifolds of the RDS.

\section{Conclusions}
\label{sec:conc}

In this paper we have shown how the action integral of classical Hamiltonian mechanics very naturally fits into the methodological framework of Lagrangian descriptors. We have applied the action based Lagrangian descriptor (LD) to standard benchmark problems, the linear saddle and the linear oscillator, for which explicit calculations can be carried out. In these examples the action based LD accurately reveals the stable and unstable manifolds of the hyperbolic equilibrium of the saddle as well as the frequency of the linear oscillator. We also considered an example of a two degree-of-freedom Hamiltonian system that arises in chemical reaction dynamics. The action based LD successfully revealed the phase spaces structures that governed the reactions dynamics, i.e. an unstable periodic orbit and its stable and unstable manifolds. We also showed how it could be used to discover Kolmogorov-Arnold-Moser tori (i.e quasiperiodic orbits) in this two degree-of-freedom Hamiltonian. In our final example, we showed that the action based LD can be used in stochastic systems. In particular, we showed how they could be used to reveal a stochastic homoclinic tangle in the stochastically forced Duffing oscillator.  It is reasonable that the connection between the action integral and Lagrangian descriptors could lead to new insights into global Hamiltonian dynamics. For example, recently it has been shown that Lagrangian descriptors can be used to understand the ‘’structure’’ of uncertainty in terms of phase space structures \cite{garcia2022bridge}. Moreover, recently ideas of symplectic topology have been used in the topic of uncertainty quantification \cite{maruskin2009dynamics, scheeres2012applications}. It would be interesting to explore the connections between these approaches more deeply.

\section*{Acknowledgments}

The authors would like to acknowledge the financial support provided by the EPSRC Grant No. EP/P021123/1 and the Office of Naval Research Grant No. N00014-01-1-0769.

\bibliography{LDs}

\end{document}